\documentclass[a4paper]{article}
\usepackage{arxiv}
\usepackage[utf8]{inputenc} 
\usepackage[T1]{fontenc}    
\usepackage{lmodern}		
\usepackage{microtype}	    
\usepackage{amsmath}
\usepackage{hyperref}    

\usepackage{url}            
\usepackage{booktabs}       
\usepackage{comment}
\newcommand{\nablax}{\nabla_{\mathbf{x}}}

\newcommand{\Deltax}{\Delta_{\mathbf{x}}}
\usepackage{amsmath}
\usepackage{amssymb}
\usepackage{amsthm}

\usepackage{nicefrac}       
\usepackage{microtype}      
\usepackage{xcolor}
\usepackage{graphicx}
\usepackage{bm}             
\usepackage[normalem]{ulem} 
\usepackage{stmaryrd}       

\usepackage{algorithmicx}
\usepackage{algorithm}
\usepackage{algpseudocode}
\usepackage{caption}
\usepackage{subcaption}
\usepackage[frozencache,cachedir=.]{minted} 

\usepackage{lineno}                   

\usepackage{multirow}
\usepackage{enumitem}
\usepackage{orcidlink}

\hypersetup{
	colorlinks=true,
	linkcolor=cyan,
	urlcolor=blue,
	citecolor=magenta,
	linktoc=all
}

\newtheorem{theorem}{Theorem}[section]
\newtheorem{proposition}[theorem]{Proposition}

\newtheorem{definition}{Definition}[section]

\newtheorem{remark}[theorem]{Remark}

\newtheoremstyle{problemstyle}
{12pt}{12pt}{\itshape}{}{\bfseries}{.}{0.5em}{}

\theoremstyle{problemstyle}
\newtheorem{tprob}{Test Problem}

\newcommand{\testproblem}[1]{%
    \begin{tprob}
        \textbf{(#1)} 
    \end{tprob}
}

\newcommand{\con}{\mathbf{u}}

\usepackage[colorinlistoftodos,textsize=tiny]{todonotes}

\title{Riemann invariant-based alternative WENO scheme for a  two-layer thin film model}
\author{
	\textbf{Biswarup Biswas} \orcidlink{0000-0001-7771-6400}\\
	Department of Mathematics\\
École Centrale School of Engineering\\
Mahindra University\\
Hyderabad, 500043, Telangana, India.\\
	\texttt{biswarupb7@gmail.com} \\
    \And
    \textbf{Rahul Barthwal} \orcidlink{0000-0002-5245-072X}\\
    Institute of Applied Analysis and Numerical Simulation,\\
     University of Stuttgart, Stuttgart, Germany\\
  \texttt{rahul.barthwal@mathematik.uni-stuttgart.de}\\
	\And
	\textbf{Rakesh Kumar} \orcidlink{0000-0002-5829-0384}\\
	Department of Mathematics\\
École Centrale School of Engineering\\
Mahindra University\\
Hyderabad, 500043, Telangana, India.\\
\texttt{rakesh.kumar@mahindrauniversity.edu.in} 
}

\makeatletter
\def\@date{}
\def\date#1{}
\makeatother
\begin{document}
\maketitle

\begin{abstract}
    In this article, we develop a multi-dimensional two-layer thin film model extending the thin film model proposed in \cite{barthwal2025hyperbolic}. The model considered in \cite{barthwal2025hyperbolic} considered a very specific Marangoni scale by choosing Marangoni numbers in both layers to be $1$. We relax this condition here and prove that the obtained system possesses a full set of Riemann invariants. Based on these findings, we develop a Riemann Invariant-based Local Characteristic Decomposition WENO (RI-WENO) method for the two-layer thin film model in one and two dimensions. The method is built upon a specially designed variable transformation constructed from the derived Riemann invariants of the system. This transformation partially diagonalizes the governing equations and yields a sparse structure in the transformed eigenvector matrices. As a result, the proposed RI-WENO framework significantly reduces the computational cost of the standard Local Characteristic Decomposition WENO approach while retaining its strong capability to suppress spurious oscillations. Numerical experiments, including new benchmark test cases, demonstrate that the RI-WENO method achieves an effective balance between accuracy and computational efficiency, making it a promising and practical choice for solving the two-layer thin film model.
\end{abstract}

\section{Introduction}
The study of thin-film equations constitutes a rich and intricate research area, especially when the flow dynamics are affected by soluble or insoluble substances that induce variations in the surface and interfacial tension. These flows appear in many technical and environmental applications, including coating technologies, thin-film solar cells, and surfactant replacement therapies (see, e.g., \cite{craster2009dynamics,  matar2004rupture, o2002theory}). Surface-tension-driven flow, usually referred to as Marangoni flow, has been widely studied over the last few decades; see, for example, \cite{matar2004rupture, jensen1992insoluble, myers1998thin} and references cited therein. Surface and interfacial tension variations may arise due to temperature gradients or due to the presence of solute particles such as surfactants and anti-surfactants. 

Mathematically, thin-film equations are typically formulated as fourth-order degenerate hyperbolic--parabolic equations, with degeneracy occurring as the film height vanishes. Even in the scalar case, the analysis of these equations and the design of robust, accurate, and high-order numerical schemes for the full lubrication models remain highly challenging. In suitable asymptotic regimes, however, the large-scale lubrication dynamics can be reduced to a first-order hyperbolic subsystem; see, e.g., \cite{barthwal2025hyperbolic, levy2006motion, cook2008shock, barthwal2025existence}. From this perspective, thin-film flows have also been extensively studied within the framework of hyperbolic balance laws, particularly in the context of gravity-driven films, inclined-plane flows, and dynamics induced by surfactants or solute particles; see, e.g., \cite{barthwal2025hyperbolic,levy2006motion, cook2008shock, barthwal2022two, bertozzi1999undercompressive, barthwal2023construction, pandey2025construction}. Nevertheless, stable and physically consistent numerical discretizations for the underlying hyperbolic subsystems remain relatively scarce, even though they seem crucial for the reliable numerical approximation of the full model.

Recently, Barthwal \& Rohde \cite{barthwal2025hyperbolic} proposed a two-layer thin film flow lubrication model, which governs the dynamics of two immiscible fluids under the influence of surface and interfacial tension induced by the presence of anti-surfactant particles. The full model turns out to be a coupled fourth-order system. However, they were able to obtain a first-order hyperbolic system from this lubrication model, which possesses a complete entropy structure and for which Riemann solutions can be constructed. Their first order system takes the following form.
\begin{equation}\label{eq: Main_system_old}
\begin{aligned}
    \dfrac{\partial f}{\partial t}+\dfrac{1}{2}\dfrac{\partial}{\partial x}\left(f^2b\right)&=0,\vspace{0.2 cm} \\
    \dfrac{\partial b}{\partial t}+\dfrac{1}{2}\dfrac{\partial}{\partial x}\left(fb^2\right)&=0,\vspace{0.2 cm} \\
     \dfrac{\partial g}{\partial t}+\dfrac{\partial}{\partial x}\left(\dfrac{g^2q}{2}+fgb\right)&=0,\vspace{0.2 cm} \\
  \dfrac{\partial q}{\partial t}+\dfrac{\partial}{\partial x}\left(\dfrac{gq^2}{2}+fbq\right)&=0,
  \end{aligned}
\end{equation}
where $f$ and $g$ denote the positive film thicknesses in the two layers.  The quantities 
$b$ 
 and
$q$ 
denote the spatial derivatives of the concentrations of the solute in each layer. 

It is noteworthy that the system \eqref{eq: Main_system_old} describes the motion of two-dimensional thin film flow under a very specific fluid regime. To be precise, Barthwal \& Rohde \cite{barthwal2025hyperbolic} considered their Marangoni numbers in both layers to be equal and to be $1$. Under this simplification, they were able to obtain explicit Riemann invariants for the system \eqref{eq: Main_system_old}, which converts the system \eqref{eq: Main_system_old} into a completely diagonal form. Even though this is a physically attainable regime, it is still very restrictive and describes only a very limited physical dynamics of the underlying system. In this article, we extend the model developed in \cite{barthwal2025hyperbolic} to its multi-dimensional counterpart while relaxing the condition on the Marangoni numbers. Precisely, the multi-dimensional first-order two-layer thin film system takes the following form:
\begin{equation}\label{eq: Main_system_new_reduced_main}
\begin{aligned}
    \dfrac{\partial f}{\partial t}+ \nablax\cdot\left(\dfrac{f^2}{2}\left(\mathrm{Ma}_1\mathbf{b}\right)\right)&=0,\vspace{0.2 cm} \\
    \dfrac{\partial \mathbf{b}}{\partial t}+\nablax\cdot \left(\dfrac{f}{2} \left(\mathrm{Ma}_1|\mathbf{b}|^2 \right)\mathbf{I}\right)&=0\\
    \dfrac{\partial g}{\partial t}+\nablax\cdot \left(fg(\mathrm{Ma}_1\mathbf{b})+\dfrac{\mathrm{Ma}_2 g^2}{2}\mathbf{q}\right)&=0\\[0.2em]
    \dfrac{\partial \mathbf{q}}{\partial t}+\nablax\cdot \left(\left(\dfrac{g\mathrm{Ma}_2 \lvert\mathbf{q}\rvert^{2}}{2}+f\left(\mathrm{Ma}_1\mathbf{b}\cdot \mathbf{q}\right)\right)\mathbf{I}\right)&=0,
 \end{aligned}
\end{equation}
where $f$ and $g$ denote the film thicknesses in the two layers.  The quantities 
$\mathbf{b}$ 
 and
$\mathbf{q}$
denote the concentration gradient vectors of the solute in each layer, and $\mathrm{Ma}_1$ and $\mathrm{Ma}_2$ are the Marangoni numbers inducing the interfacial and surface tension, respectively. To study the hyperbolic structure of the model, we analyze the system \eqref{eq: Main_system_new_reduced_main} by projecting the multidimensional equations onto a one-dimensional direction. This allows us to explicitly characterize the hyperbolicity of the resulting projected system. Moreover, the mathematical structure of \eqref{eq: Main_system_new_reduced_main} makes it particularly well suited for the construction of higher-order accurate numerical schemes. In fact, when the tangential component of the concentration gradient vanishes, we are able to derive a complete set of explicit Riemann invariants for \eqref{eq: Main_system_new_reduced_main}. These invariants diagonalize the system completely and therefore provide a natural framework for the design of higher-order numerical methods based on reconstruction in Riemann invariant variables.

Among high-order methods for problems with discontinuities, Weighted Essentially Non-Oscillatory (WENO) schemes have become one of the most widely used approaches, particularly for hyperbolic conservation laws. The WENO methodology was first introduced in the finite volume framework in \cite{liu1994weighted} and later extended to the finite difference setting in \cite{jiang1996efficient}. Over the past three decades, WENO schemes have undergone significant development, both theoretically and in their application to increasingly complex problems. Several improved variants have been proposed, including WENO-JS \cite{jiang1996efficient}, WENO-Z \cite{ borges2008improved}, and WENO-AO \cite{bal-etal_16a, kum-cha_18a, kum-cha_19a}, among others. Beyond the classical finite difference and finite volume formulations, WENO methodologies have also been extended to alternative frameworks that incorporate Riemann solvers within a finite difference setting \cite{balsara2025efficient}.
While WENO schemes perform remarkably well for scalar problems, their extension to systems of equations may lead to spurious oscillations \cite{jiang1996efficient,kum_cha_22a,aru_etal_22a,MR4946687}. A common strategy to mitigate these oscillations is to perform reconstruction in the characteristic space \cite{jiang1996efficient}. However, this approach necessitates the computation of eigenvectors, thereby increasing the computational cost.
To alleviate this issue, several hybrid strategies have been developed, wherein characteristic-wise reconstruction is employed near discontinuities, while component-wise reconstruction is used in smooth regions to reduce computational overhead \cite{kum_cha_22a, aru_etal_22a, zhao-etal_20a}. More recently, methods based on local characteristic decomposition using Riemann invariants have been proposed, demonstrating significant computational savings without compromising accuracy \cite{xu2024local,wu2025finite}.

Motivated by these developments, we exploit the Riemann-invariant structure of the governing system \eqref{eq: Main_system_new_reduced_main} to design an alternative WENO reconstruction in Riemann-invariant variables. Although the system admits a set of Riemann invariants that forms a coordinate system, recovering the physical variables through the inverse transformation requires solving a nonlinear quartic equation, which increases the computational cost of the scheme. To avoid this difficulty, we instead introduce a modified set of Riemann invariants that yields a sparser characteristic decomposition. This choice not only reduces the computational expense but also improves the accuracy of the numerical method. Moreover, we introduce new benchmark test cases in both one and two space dimensions for the system \eqref{eq: Main_system_new_reduced_main}. In particular, we investigate the effect of two distinct Marangoni numbers on the flow. Moreover, under the assumption that the corresponding concentration gradient vectors are aligned with the line $y=x$, we are able to employ the Riemann invariant structure even in the two-dimensional setting. Physically, this assumption implies that the concentration varies at the same rate in the $x$- and $y$-directions, so that no coordinate direction is preferred. Consequently, the steepest variation occurs along the diagonal direction, reflecting a symmetric influence of the concentration gradient in both spatial directions.

The rest of the article is structured as follows. In section 2, we formulate the multi-dimensional two-layer thin film model and obtain the governing system \eqref{eq: Main_system_new_reduced_main}. Section 3 is devoted to discuss the hyperbolicity and the underlying Riemann invariant structure of the system \eqref{eq: Main_system_new_reduced_main}. In Section 4, we define the numerical schemes and complement them with some numerical examples in Section 5. Conclusions and future outlook are provided in Section 6.
\section{Multi-dimensional two-layer thin film model under the influence of a perfectly soluble anti-surfactant\label{sec:2}}
In this section, we develop the multi-dimensional version of the hyperbolic model developed in \cite{barthwal2025hyperbolic}, which governs the first-order dynamics of a two-layer thin film flow under the influence of a perfectly soluble anti-surfactant solute. As the model considered here is a natural extension of the model developed in \cite{barthwal2025hyperbolic}, we leave most of the physical details of the model and refer the interested reader to \cite{barthwal2025hyperbolic, barthwal2026generalized}; see also \cite{barthwal2025existence, barthwal2023construction, conn2016fluid}.
\subsection{Evolution equations for the film heights and concentrations}
Here, we consider the model considered in \cite{barthwal2023construction} and \cite{conn2016fluid}, but consider the geometry such that the substrate is on a two-dimensional plane. The free surface is assumed to be located at a height $z=(f+g)(x, y, t)$, while the interface is assumed to be located at $z=g(x, y, t)$. This, in particular, implies that the model here governs the three-dimensional thin film flows. Following \cite{barthwal2025hyperbolic, barthwal2023construction}, we assume that $\epsilon=H/L<\!<\! 1$ for thin film flow, where $H$ and $L$ are characteristic height and length, respectively. We use a similar non-dimensional scaling as in \cite{barthwal2025hyperbolic, barthwal2023construction} and \cite{matar2004rupture}, and neglect the terms of order $\epsilon^2$ in the two-phase incompressible Navier-Stokes equations, along with the concentration transport equations for a perfectly soluble solute. This then leads to a reduced system of governing equations of the form
\begin{subequations}
\begin{align}
    \nabla_{\mathbf{x}}\cdot  \mathbf{u}_i+\partial_z w_i&=0, \label{eq: 2.3a} \vspace{0.2 cm}\\
   \nabla_{\mathbf{x}} p_i&=\partial_z^2\,\mathbf{u}_i,\label{eq: 2.3b} \vspace{0.2 cm}\\
  \partial_z p_i&=0, \label{eq: 2.3d} \vspace{0.2 cm}\\
   \partial_t c_i^0+\mathbf{u}_i\cdot\nabla c_i^0&=(\mathrm{Pb}_i)^{-1}\left(\Delta_{\mathbf{x}} c_i^0+\partial^2_z c_i^1\right).\label{eq: 2.3e} 
\end{align}
\end{subequations}
In \eqref{eq: 2.3a}-\eqref{eq: 2.3e}, $\mathbf{u}_i=(u_i, v_i)$ is the $x-y$ components while $w_i$ is the $z$-component of the non-dimensional velocity vectors, $p_i$ is the pressure in the $ith$ layer for $i\in \{1, 2\}$. $\mathbf{x}=(x, y)$ is the horizontal coordinates such that $\nabla_{\mathbf{x}}=(\partial_x, \partial_y)$ denotes the gradient operator while $\Delta_{\mathbf{x}}=(\partial_x^2, \partial_y^2)$ is the Laplacian operator. Moreover, $\mathrm{Pb}_i\, , i\in\{1, 2\}$ are the P\'eclet numbers with $\mathrm{Pb}_i=O(1)$ and the concentrations $c_i$ are assumed to satisfy an expansion of the form \cite{barthwal2025hyperbolic, barthwal2023construction}
\[
c_i=c_i^0(\mathbf{x}, t)+\epsilon^2 c_i^1(\mathbf{x}, z, t)+O(\epsilon^4).
\]
Neglecting the terms of order $\epsilon^2$ and simplifying, we obtain non-dimensional boundary conditions of the form \cite{barthwal2025hyperbolic, barthwal2023construction}
\begin{subequations}\label{eq: 2.17main}
\begin{alignat}{2}
(\mathbf{u_1},w_1)&=(\mathbf{0}, 0)        &\quad \mathrm{on}~~ \{z=0\},\label{eq: 2.4a}\\
\mathbf{u_1}&=\mathbf{u_2}       &\quad \mathrm{on}~~ \{z=f\}, \label{eq: 2.4b}\vspace{0.1 cm}\\
\partial_t f+\mathbf{u_1}\cdot \nabla f&=w_1        &\quad \mathrm{on}~~ \{z=f\}, \label{eq: 2.4c}\vspace{0.1 cm}\\
\partial_t g+\mathbf{u_2}\cdot \nabla g&=w_2        &\quad\quad \mathrm{on}~~ \{z=f+g\}, \label{eq: 2.4d}\vspace{0.1 cm}\\
\mu p_2-p_1&=(\mathrm{Ca_1})^{-1}\Delta_{\mathbf{x}}^2 f       &\quad \mathrm{on}~~\{z=f\},\label{eq: 2.4e} \vspace{0.1 cm}\\
-p_2&=(\mathrm{Ca_2})^{-1}\Delta_{\mathbf{x}}^2 (f+g)       &\quad \mathrm{on}~~\{z=f+g\},\label{eq: 2.4f} \vspace{0.1 cm}\\
\partial_z \mathbf{u}_1&=\mu\partial_z \mathbf{u}_2+\mathrm{Ma_1} \nablax c_1       &\quad \mathrm{on}~~ \{z=f\},\label{eq: 2.4g}\\
\partial_z \mathbf{u}_2&=\mathrm{Ma_2}\nablax c_2        &\quad \mathrm{on}~ \{z=f+g\},\label{eq: 2.4h} \\
\partial_z c_1^1
    &= 0,                                      &\quad \mathrm{on}~\{z=0\},       \label{eq: 2.4i}\\
\partial_z c_1^1
    &= \nabla_{\mathbf{x}}c_1^0\cdot \nabla_{\mathbf{x}} f,        &\quad \mathrm{on}~\{z=f\},       \label{eq: 2.4j}\\
\partial_z c_2^1
    &= \nabla_{\mathbf{x}}c_2^0\cdot \nabla_{\mathbf{x}} f,        &\quad \mathrm{on}~\{z=f\},       \label{eq: 2.4k1}\\
\partial_z c_2^1
    &= \nabla_{\mathbf{x}}c_2^0\cdot \nabla_{\mathbf{x}} (f+g),          &\quad \mathrm{on}~\{z=f+g\}.     \label{eq: 2.4k2}
 \end{alignat}
\end{subequations}
In \eqref{eq: 2.17main}, $\mu=\mu_2^*/\mu_1^*$ denotes the ratio of the viscosities of two fluids, $\mathrm{Ma_i}$ denote Marangoni numbers, and $\mathrm{Ca_i}$ are the capillarity numbers all of order 1 (see \cite{barthwal2025hyperbolic} for more details on the scalings of these dimensionless numbers).
\subsection{Evolution equation for the film thicknesses}
Using the kinematic boundary conditions \eqref{eq: 2.4c}-\eqref{eq: 2.4d}, the evolution equation of film thicknesses $f$ and $g$ are then given by
\begin{align}
    \dfrac{\partial f}{\partial t}+\left(\displaystyle\int_{0}^f u_1 dz\right)+\left(\displaystyle\int_{0}^f v_1 dz\right)=0,\label{evolution_film_height_first}\\
    \dfrac{\partial g}{\partial t}+\left(\displaystyle\int_{f}^{f+g} u_2 dz\right)+\left(\displaystyle\int_{f}^{f+g} v_2 dz\right)=0.\label{evolution_film_height_second}
\end{align}
Now using \eqref{eq: 2.3d} and normal stress balances \eqref{eq: 2.4e}-\eqref{eq: 2.4f}, we have
\begin{align*}
p_2(\mathbf{x}, z, t)=-(\mathrm{Ca_2})^{-1}\Delta_{\mathbf{x}}^2 (f+g),~ \text{and}\quad p_1(\mathbf{x}, z, t)=-\mu (\mathrm{Ca_2})^{-1}\Delta_{\mathbf{x}}^2 (f+g)-(\mathrm{Ca_1})^{-1}\Delta_{\mathbf{x}}^2 f \quad \forall (\mathbf{x}. z)\in \Omega
\end{align*}
Hence, by \eqref{eq: 2.3b}, we have
\begin{align*}
\partial_z \mathbf{u}_1=\nabla_x p_1\,z+a_1. 
\end{align*}
Thus, using the no-slip condition \eqref{eq: 2.4a}, we get
\begin{align}\label{eq: 2.22mainb}
   \mathbf{u}_1=\nabla_x p_1\dfrac{z^2}{2}+a_1z. 
\end{align}
Now from the tangential stress balance \eqref{eq: 2.4g}, we obtain
\begin{align}\label{a_1}
    a_1= \mu \partial_z \mathbf{u}_2\bigg|_{z=f}+\mathrm{Ma_1} \nablax c_1-\bigg[\nablax p_1\bigg]f.
\end{align}
Moreover, again from \eqref{eq: 2.3b}, we have 
\[
\partial_z \mathbf{u}_2=z\nablax p_2+a_2.
\]
Thus, using the tangential stress balance we get
$$a_2= \mathrm{Ma_2}\nablax c_2-\nablax p_2(f+g)$$
or 
\begin{align}\label{eq: 2.20u_2}
\partial_z \mathbf{u}_2=-\nablax p_2((f+g)-z)+\mathrm{Ma_2}\nablax c_2,
\end{align}
which implies 
\begin{align*}
    \partial_z \mathbf{u}_2\bigg|_{z=f}=-g\,\nablax p_2 +\mathrm{Ma_2}\nablax c_2. 
\end{align*}
Hence by \eqref{a_1}, we have
\begin{align*}
    a_1&=-\mu g\,\nablax p_2 +\mu \mathrm{Ma_2}\nablax c_2+\mathrm{Ma_1}\nablax c_1-f\,\nablax p_1.
\end{align*}
Therefore, by \eqref{eq: 2.22mainb}, we obtain
\begin{align}\label{eq: u_1}
   \mathbf{u}_1=-\mu \nablax p_2 \,gz+\mu \mathrm{Ma_2}\nablax c_2 \,z-\nablax p_1\,z\left(f-\dfrac{z}{2}\right)+\mathrm{Ma_1} \nabla_x c_1\,z.
\end{align}
Now \eqref{eq: 2.20u_2} implies
\begin{align*}
\mathbf{u}_2&=-\nablax p_2\,z\left((f+g)-\dfrac{z}{2}\right)+\mathrm{Ma_2}\nablax c_2z+a_3.
\end{align*}
At $\{z=f\}$, we have $\mathbf{u}_2=\mathbf{u}_1$ from \eqref{eq: 2.4a}, therefore
\begin{align*}
a_3=-\nablax p_2\,f\left(g(\mu-1)-\dfrac{f}{2}\right)-\nablax p_1\dfrac{f^2}{2}+\mathrm{Ma_2}\nablax c_2\,f(\mu-1)+\mathrm{Ma_1} \nablax c_1\,f,
\end{align*}
and hence
\begin{align}\label{eq:u_2}
\mathbf{u}_2&=-\nablax p_2\,\bigg[f\left(g(\mu-1)-\dfrac{f}{2}\right)+z\left((f+g)-\dfrac{z}{2}\right)\bigg]-\nablax p_1\dfrac{f^2}{2} \nonumber\\
&\hspace{3.5 cm}+\mathrm{Ma_2}\nablax c_2\left(f(\mu-1)+z\right)+\mathrm{Ma_1} \nablax c_1\,f.
\end{align}
Then, by integrating \eqref{eq: u_1} from $z=0$ to $z=f$ and \eqref{eq: 2.20u_2} from $z=f$ to $z=f+g$ and using \eqref{evolution_film_height_first} and \eqref{evolution_film_height_second}, we obtain the evolution equations for the film heights as
\begin{align}\label{eq: film_height}
    \dfrac{\partial f}{\partial t}+\nablax\cdot \left(\dfrac{ f^2}{2}\nablax(\mathrm{Ma}_1c_1+ \mathrm{Ma}_2\mu c_2)\right)=\nablax\cdot \mathbf{Q}_1,\\
     \dfrac{\partial g}{\partial t}+\nablax\cdot \left(fg\nablax(\mathrm{Ma}_1c_1+ \mathrm{Ma}_2\mu c_2)+\dfrac{\mathrm{Ma}_2 g^2}{2}\nablax c_2\right)=\nablax\cdot \mathbf{Q}_2,
\end{align}
where
\begin{align}\label{eq: 2.22main}
    &\mathbf{Q}_1=-\dfrac{\mu f^2\nablax \Deltax(f+g)}{2\mathrm{Ca_2}}\left(\dfrac{2f}{3}+g\right)-\bigg[\dfrac{f^3\nablax \Deltax f}{3\mathrm{Ca_1}}\bigg]\\
&\mathbf{Q}_2=-\dfrac{\nablax \Deltax(f+g)}{\mathrm{Ca_2}}\bigg[\left(\dfrac{(\mu-1)(f^2 g+2fg^2)}{2}\right)+\dfrac{(f+g)^3}{3}-\dfrac{f^2(2f+3g)}{6}\bigg]-\dfrac{f^2 g~ \nablax \Deltax f}{2~\mathrm{Ca_1}}\label{eq: 2.23main}.
\end{align}
\subsection{Evolution equations for concentrations in each layer}
The evolution equation for the leading order term $c_i^0$ can be obtained by integrating \eqref{eq: 2.3e} with respect to $z$ from $\{z=0\}$ to $\{z=f\}$ for $i=1$ and from $\{z=f\}$ to $\{z=f+g\}$ for $i=2$ with boundary conditions \eqref{eq: 2.4i}-\eqref{eq: 2.4k2}. The evolution equation of $c_i^0(\mathbf{x}, t), i\in \{1, 2\}$ reads as
\begin{align}\label{eq: concentration}
    \partial_t c_1^0+\left(\dfrac{ f}{2}\nablax(\mathrm{Ma}_1c_1^0+ \mu\mathrm{Ma}_2 c_2^0)\right)\cdot \nablax c_1^0=\dfrac{(\mathrm{Pb}_1)^{-1}}{f}\nabla \cdot(f\nablax c_1^0)+\mathbf{Q}_1\cdot \nablax c_1^0,\\
    \partial_t c_2^0+\left(f\nablax(\mathrm{Ma}_1c_1^0+ \mu\mathrm{Ma}_2 c_2^0)+\dfrac{\mathrm{Ma}_2 g}{2}\nablax c_2^0\right)\cdot\nablax c_2^0=\dfrac{(\mathrm{Pb}_2)^{-1}}{g}\nabla \cdot(g\nablax c_2^0)+\mathbf{Q}_2\cdot \nablax c_2^0,\label{eq: concentration_2}
\end{align}
where $\mathbf{Q}_i$, $i\in \{1, 2\}$ are the higher order nonlinear terms as defined in \eqref{eq: 2.22main}-\eqref{eq: 2.23main}.
\subsection{Reduced evolution equations for first-order dynamics}
Following \cite{barthwal2025hyperbolic, barthwal2026generalized, conn2017simple}, we assume that the capillarity and diffusivity effects are negligible. Under these simplified assumptions, we can leave the higher-order terms from the evolution equations \eqref{eq: film_height}-\eqref{eq: concentration_2} and obtain reduced first-order governing equations of the form
\begin{equation}\label{eq:2.9new}
\begin{aligned}
 \frac{\partial f}{\partial t}
 +\nablax\cdot \left(\dfrac{ f^2}{2}\nablax(\mathrm{Ma}_1c_1+ \mu \mathrm{Ma}_2c_2)\right) &= 0,\\[0.2em]
  \dfrac{\partial g}{\partial t}+\nablax\cdot \left(fg\nablax(\mathrm{Ma}_1c_1+ \mu\mathrm{Ma}_2 c_2)+\dfrac{\mathrm{Ma}_2 g^2}{2}\nablax c_2\right)&=0\\[0.2em]
 \frac{\partial c_1}{\partial t}
 + \dfrac{f}{2}\left(\mathrm{Ma}_1\lvert\nablax c_1\rvert^{2}+\mu\mathrm{Ma}_2 \nablax c_1\cdot \nablax c_2\right) &= 0,\\[0.2em]
 \frac{\partial c_2}{\partial t}
 + \dfrac{g\mathrm{Ma}_2 \lvert\nablax c_2\rvert^{2}}{2}+f\left(\mathrm{Ma}_1\lvert\nablax c_1\rvert^{2}+\mu\mathrm{Ma}_2 \nablax c_1\cdot \nablax c_2\right) &= 0,
\end{aligned}
\end{equation}
where we have dropped the subscript from $c_i^0$ for $i\in \{1, 2\}$.
\subsection{A purely conservative first order subsystem}
To uncover the conservative structure of the reduced evolution equations \eqref{eq:2.9new}, we take gradient of the concentration evolution equations, which gives us for the film thicknesses $f$, $g$ and the concentration gradient $\mathbf{b}=\nablax c_1=\left(\partial_x c_1, \partial_y c_1\right)^\top=(b_1, b_2)^\top$ and $\mathbf{q}=\nablax c_2=\left(\partial_x c_2, \partial_y c_2\right)^\top=(q_1, q_2)^\top$
the following set of nonlinear, first-order hyperbolic equations. 
\begin{equation}\label{eq: Main_system_new_updated}
\begin{aligned}
    \dfrac{\partial f}{\partial t}+ \nablax\cdot\left(\dfrac{f^2}{2}\left(\mathrm{Ma}_1\mathbf{b}+\mu \mathrm{Ma}_2\mathbf{q}\right)\right)&=0,\vspace{0.2 cm} \\
    \dfrac{\partial \mathbf{b}}{\partial t}+\nablax\cdot \left(\dfrac{f}{2} \left(\mathrm{Ma}_1|\mathbf{b}|^2 +\mu \mathrm{Ma}_2\mathbf{b}\cdot \mathbf{q})\right)\mathbf{I}\right)&=0\\
    \dfrac{\partial g}{\partial t}+\nablax\cdot \left(fg(\mathrm{Ma}_1\mathbf{b}+ \mu\mathrm{Ma}_2 \mathbf{q})+\dfrac{\mathrm{Ma}_2 g^2}{2}\mathbf{q}\right)&=0\\[0.2em]
    \dfrac{\partial \mathbf{q}}{\partial t}+\nablax\cdot \left(\left(\dfrac{g\mathrm{Ma}_2 \lvert\mathbf{q}\rvert^{2}}{2}+f\left(\mathrm{Ma}_1\mathbf{b}\cdot \mathbf{q}+\mu\mathrm{Ma}_2 \lvert\mathbf{q}\rvert^2\right)\right)\mathbf{I}\right)&=0.
 \end{aligned}
\end{equation}
Without loss of generality, we take $\mathrm{Ma}_1=\alpha>0$ and $\mathrm{Ma}_2=\beta>0$ . Moreover, we choose the viscosity ratio very low, which is typically true in oil-water type flows, such that $\mu=0$. With these simplifications, the system reduces to 
\begin{equation}\label{eq: Main_system_new_reduced}
\begin{aligned}
    \dfrac{\partial f}{\partial t}+ \nablax\cdot\left(\dfrac{f^2}{2}\left(\alpha\mathbf{b}\right)\right)&=0,\vspace{0.2 cm} \\
    \dfrac{\partial \mathbf{b}}{\partial t}+\nablax\cdot \left(\dfrac{f}{2} \left(\alpha|\mathbf{b}|^2 \right)\mathbf{I}\right)&=0\\
    \dfrac{\partial g}{\partial t}+\nablax\cdot \left(fg(\alpha\mathbf{b})+\dfrac{\beta g^2\mathbf{q}}{2}\right)&=0\\[0.2em]
    \dfrac{\partial \mathbf{q}}{\partial t}+\nablax\cdot \left(\left(\dfrac{\beta g \lvert\mathbf{q}\rvert^{2}}{2}+f\left(\alpha\mathbf{b}\cdot \mathbf{q}\right)\right)\mathbf{I}\right)&=0.
 \end{aligned}
\end{equation}
\begin{remark}
It is noteworthy that the system \eqref{1d_system} reduces to the one-dimensional system considered in \cite{barthwal2025hyperbolic} for $\alpha=\beta=1$. This indicates the consistency of our modelling approach. 
\end{remark}
\section{Preliminaries, Hyperbolicity of the System~\eqref{eq: Main_system_new_reduced},
and Riemann Invariants}%
\label{sec: hyperbolicity}
 
\subsection{Preliminaries}
 
In this section, we recall some foundational facts on hyperbolicity and
Riemann invariants for systems of conservation laws, following
\cite{wu2025finite, dafermos2005hyperbolic, serre1999systems}.
 
Consider a system of conservation laws of $m$ components in one space dimension,
\begin{equation}\label{hyp_1}
  \mathbf{u}_t + \mathbf{f}(\mathbf{u})_x = \mathbf{0},
\end{equation}
where $\mathbf{u} = (u_1, u_2, \ldots, u_m)^T \in \mathbb{R}^m$ are the
conserved variables and
$\mathbf{f}(\mathbf{u}) = (f_1(\mathbf{u}), f_2(\mathbf{u}), \ldots, f_m(\mathbf{u}))^T
\in C^1(\mathbb{R}^m)$ are the flux functions.
For smooth solutions, the system~\eqref{hyp_1} can be rewritten in the
quasilinear form
\begin{equation}\label{hyp_2}
  \mathbf{u}_t + \mathbf{A}(\mathbf{u})\,\mathbf{u}_x = \mathbf{0},
\end{equation}
where $\mathbf{A}(\mathbf{u}) = \dfrac{\partial \mathbf{f}}{\partial \mathbf{u}}$
is the Jacobian matrix.
 
\begin{definition}[Hyperbolicity in 1D]\label{hyper_1_d}
The quasilinear system~\eqref{hyp_2} is \emph{hyperbolic} if the Jacobian
matrix $\mathbf{A}(\mathbf{u})$ is diagonalizable for all $\mathbf{u}$, that is, there exists an invertible matrix $\mathbf{R}(\mathbf{u})$ such that
\[
  \mathbf{A}(\mathbf{u})
  = \mathbf{R}(\mathbf{u})\,\boldsymbol{\Lambda}(\mathbf{u})\,
    \mathbf{R}(\mathbf{u})^{-1},
\]
where
\[
  \boldsymbol{\Lambda}(\mathbf{u})
  = \operatorname{diag}\!\bigl(\lambda_1(\mathbf{u}),\,\lambda_2(\mathbf{u}),
    \ldots,\,\lambda_m(\mathbf{u})\bigr)
\]
is the diagonal matrix of real eigenvalues,
$\mathbf{R}(\mathbf{u}) = \bigl(\mathbf{r}_1(\mathbf{u}),\,
\mathbf{r}_2(\mathbf{u}),\,\ldots,\,\mathbf{r}_m(\mathbf{u})\bigr)$
is the right eigenmatrix whose columns are the right eigenvectors, and
$\mathbf{R}(\mathbf{u})^{-1} = \mathbf{L}(\mathbf{u}) =
\bigl(\mathbf{l}_1(\mathbf{u})^T,\,\mathbf{l}_2(\mathbf{u})^T,\,\ldots,\,
\mathbf{l}_m(\mathbf{u})^T\bigr)^T$
is the left eigenmatrix whose rows are the left eigenvectors.
\end{definition}
 
The one-dimensional notion extends naturally to multiple spatial dimensions.
 
\begin{definition}[Hyperbolicity in multiple dimensions]\label{hyper_multi_d}
A system of conservation laws in $d$ spatial dimensions,
\[
  \mathbf{u}_t + \sum_{i=1}^{d} \mathbf{f}_i(\mathbf{u})_{x_i} = \mathbf{0},
\]
is called \emph{hyperbolic} if, for every unit vector
$\mathbf{n} = (n_1, n_2, \ldots, n_d)^T \in \mathbb{R}^d$, the associated
one-dimensional problem
\[
  \mathbf{u}_t + \Bigl(\sum_{i=1}^{d} n_i\,\mathbf{f}_i(\mathbf{u})\Bigr)_x
  = \mathbf{0}
\]
is hyperbolic in the sense of Definition~\ref{hyper_1_d}.
\end{definition}
 
\subsection{Change of Variables and Transformed Jacobian}
 
We next consider variable transformations derived from the eigen-structure of
the Jacobian $\mathbf{A}(\mathbf{u})$.
Such transformations play a fundamental role in the Local Characteristic
Decomposition (LCD) procedure described in Section~\ref{sec: Numerical_scheme}.
 
Let $\mathbf{V}(\mathbf{u}) : \mathcal{U} \to \mathbb{R}^m$ be a
differentiable and invertible transformation defined on the state space
$\mathcal{U} \subseteq \mathbb{R}^m$.
The system~\eqref{hyp_1} can then be rewritten in the new variables
$\mathbf{V}$ as
\begin{equation}\label{eq:transformed_system}
  \mathbf{V}_t + \tilde{\mathbf{A}}(\mathbf{V})\,\mathbf{V}_x = \mathbf{0},
\end{equation}
where the transformed Jacobian matrix is
\begin{equation}\label{eq:transformed_Jacobian}
  \tilde{\mathbf{A}}
  = \frac{\partial \mathbf{V}}{\partial \mathbf{u}}\,
    \mathbf{A}\,
    \frac{\partial \mathbf{u}}{\partial \mathbf{V}}.
\end{equation}
 
\begin{proposition}[Change of variables]\label{prop:change_of_vars}
If the transformation $\mathbf{V}(\mathbf{u})$ is differentiable and
invertible, then the transformed Jacobian $\tilde{\mathbf{A}}$ has the same
eigenvalues $\lambda_k$ as the original Jacobian $\mathbf{A}$.
Moreover, the right and left eigenmatrices transform according to
\begin{equation}
  \tilde{\mathbf{R}}
  = \frac{\partial \mathbf{V}}{\partial \mathbf{u}}\,\mathbf{R},
  \qquad
  \tilde{\mathbf{L}}
  = \mathbf{L}\,\frac{\partial \mathbf{u}}{\partial \mathbf{V}}.
\end{equation}
\end{proposition}

 
\subsection{Riemann Invariants}
 
The structure of the transformed Jacobian $\tilde{\mathbf{A}}$ depends
critically on how the new variables relate to the characteristic fields of the
system.
This relationship is made precise through the notion of a Riemann invariant.
 
\begin{definition}[Riemann invariant]\label{Defn_RI}
A scalar function $w : \mathcal{U} \to \mathbb{R}$ is called a
\emph{$k$-Riemann invariant} of the hyperbolic system~\eqref{hyp_1} if it
satisfies
\begin{equation}\label{eq:RI_def}
  \nabla_{\mathbf{u}}\,w(\mathbf{u}) \cdot \mathbf{r}_k(\mathbf{u}) = 0,
  \qquad \forall\,\mathbf{u} \in \mathcal{U},
\end{equation}
where $\mathbf{r}_k(\mathbf{u})$ is the $k$-th right eigenvector of
$\mathbf{A}(\mathbf{u})$.
\end{definition}

The following result \cite{xu2024local}, which is central to the design of characteristic-based
numerical schemes, shows that choosing Riemann invariants as the transformation
variables introduce zeros into the transformed Jacobian.
 
\begin{proposition}[Sparsity of the transformed Jacobian]
\label{prop:sparsity}
Let $\mathbf{V}(\mathbf{u}) = (v_1, v_2, \ldots, v_m)^T$ be a differentiable
and invertible transformation, and let $\tilde{\mathbf{A}}$ be the transformed
Jacobian defined in~\eqref{eq:transformed_Jacobian}.
If the component $v_i(\mathbf{u})$ is a $k$-Riemann invariant, then the
$(i,k)$-th entry of $\tilde{\mathbf{A}}$ vanishes:
\[
  \tilde{A}_{ik} = 0.
\]
\end{proposition}
 
The following remark identifies the strongest possible consequence of
Proposition~\ref{prop:sparsity}: a coordinate system composed entirely of
Riemann invariants completely diagonalizes the system.
 
\begin{remark}[Full diagonalization via Riemann invariants]
\label{rem:diagonalization}
Suppose there exists an invertible transformation
$\mathbf{V}(\mathbf{u}) = (v_1, v_2, \ldots, v_m)^T$ such that, for each
index $i \in \{1, \ldots, m\}$, the component $v_i$ is a $k$-Riemann
invariant for every $k \neq i$.
Then, by Proposition~\ref{prop:sparsity}, every off-diagonal entry of
$\tilde{\mathbf{A}}$ vanishes, and the transformed Jacobian reduces to
\[
  \tilde{\mathbf{A}}
  = \operatorname{diag}\!\bigl(\lambda_1(\mathbf{u}),\,\lambda_2(\mathbf{u}),
    \ldots,\,\lambda_m(\mathbf{u})\bigr).
\]
In this coordinate system, the quasilinear system~\eqref{eq:transformed_system}
decouples completely into $m$ independent scalar advection equations,
each governed by a single characteristic speed $\lambda_i$.
 
Such a global coordinate system of Riemann invariants always exists for
hyperbolic systems of two equations \cite{dafermos2005hyperbolic}.
For systems of three or more equations, it exists only in special cases;
in general, only local (i.e.\ in a neighbourhood of a point in state space)
coordinates of this type can be guaranteed, and only under genuine nonlinearity
or linear degeneracy conditions on the characteristic fields
\cite{dafermos2005hyperbolic, serre1999systems}.
When such a coordinate system is available, it is the natural setting for
the construction of Riemann solvers and high-order characteristic-based
schemes.
\end{remark}
\subsection{Hyperbolicity of the system \eqref{eq: Main_system_new_reduced}}
In what follows, we prove that the two-dimensional system \eqref{eq: Main_system_new_reduced} is a hyperbolic system of conservation laws in the sense of Definition \ref{hyper_multi_d}. Let $\mathbf{n}:=(n_1, n_2)\in \mathbb{R}^2$ be a unit normal vector and let $\mathbf{t}:=(t_1, t_2)\in \mathbb{R}^2$ be a unit tangent vector such that $\{\mathbf{n}, \mathbf{t}\}$ form an orthonormal basis of $\mathbb{R}^2$. Then we can decompose the concentration gradient vectors $\mathbf{b}$ and $\mathbf{q}$ as a linear combination of these vectors. Precisely, we can express
\[
\mathbf{b}=(\mathbf{b}\cdot \mathbf{n})\mathbf{n}+(\mathbf{b}\cdot \mathbf{t})\mathbf{t},\qquad \qquad \mathbf{q}=(\mathbf{q}\cdot \mathbf{n})\mathbf{n}+(\mathbf{q}\cdot \mathbf{t})\mathbf{t}
\]
Assuming that the solutions of the system \eqref{eq: Main_system_new_reduced} depend only on the direction $x=\mathbf{x}\cdot \mathbf{n}$, we have
\[
\nablax:=\mathbf{n}\partial_x.
\]
Therefore, one can rewrite the system \eqref{eq: Main_system_new_reduced} as the one-dimensional directional system
\begin{equation}\label{eq:directional_system}
\begin{aligned}
f_t+\partial_x\left(\frac{\alpha}{2}f^2b\right)&=0,\\[0.2em]
b_t+\partial_x\left(\frac{\alpha}{2}f(b^2+\eta^2)\right)&=0,\\[0.2em]
\eta_t&=0,\\[0.2em]
g_t+\partial_x\left(\alpha fb g+\frac{\beta}{2}g^2q\right)&=0,\\[0.2em]
q_t+\partial_x\left(\frac{\beta}{2}g(q^2+\zeta^2)+\alpha f(bq+\eta\zeta)\right)&=0,\\[0.2em]
\zeta_t&=0,
\end{aligned}
\end{equation}
where $b=\mathbf{b}\cdot \mathbf{n}$, $q=\mathbf{q}\cdot \mathbf{n}$, $\eta=\mathbf{b}\cdot \mathbf{t}$ and $\zeta=\mathbf{q}\cdot \mathbf{t}$.

Hence, in the variables $\mathbf{u_n}=(f,b,\eta,g,q,\zeta)^T$, the system \eqref{eq:directional_system} can be written in a compact form as
\[
\mathbf{u}_t+\mathbf{F_n(u)}_{x}=0
\]
where
\[
\mathbf{F_n(u)}=
\left(
\dfrac{\alpha}{2}f^2b,~
\dfrac{\alpha}{2}f(b^2+\eta^2),~
0,~
\alpha fgb+\dfrac{\beta}{2}g^2q,~
\dfrac{\beta}{2}g(q^2+\zeta^2)+\alpha f(bq+\eta\zeta),~
0\right)^\top.
\]
with the flux Jacobian 
\[
\mathbf{A}_{\mathbf{n}}(\mathbf{u})
=
\frac{\partial \mathbf{F}_{\mathbf{n}}(\mathbf{u})}{\partial \mathbf{u}}
=
\begin{pmatrix}
\alpha fb & \dfrac{\alpha}{2}f^2 & 0 & 0 & 0 & 0\\[0.2em]
\dfrac{\alpha}{2}(b^2+\eta^2) & \alpha fb & \alpha f\eta & 0 & 0 & 0\\[0.2em]
0 & 0 & 0 & 0 & 0 & 0\\[0.2em]
\alpha gb & \alpha fg & 0 & \alpha fb+\beta gq & \dfrac{\beta}{2}g^2 & 0\\[0.2em]
\alpha (bq+\eta\zeta) & \alpha fq & \alpha f\zeta & \dfrac{\beta}{2}(q^2+\zeta^2) & \alpha fb+\beta gq & \beta g\zeta\\[0.2em]
0 & 0 & 0 & 0 & 0 & 0
\end{pmatrix}.
\]
We note that the Jacobian $\mathbf{A}_{\mathbf{n}}(\mathbf{u})$ is a $6\times 6$ matrix, corresponding to the six variables $(f,b,\eta,g,q,\zeta).$ In particular, the tangential components $\eta$ and $\zeta$ satisfy
\[
\eta_t=0,\qquad \zeta_t=0,
\]
and therefore give rise to two zero characteristic speeds. Thus, the full directional system has six eigenvalues, which include four nontrivial characteristic speeds and two zeros.

To identify the four nonzero eigenvalues, observe that the $4\times 4$ block is
\[
\mathbf{A}_{4\times 4}=
\begin{pmatrix}
\alpha fb & \dfrac{\alpha}{2}f^2 & 0 & 0\\[0.2em]
\dfrac{\alpha}{2}(b^2+\eta^2) & \alpha fb & 0 & 0\\[0.2em]
\alpha gb & \alpha fg & \alpha fb+\beta gq & \dfrac{\beta}{2}g^2\\[0.2em]
\alpha (bq+\eta\zeta) & \alpha fq & \dfrac{\beta}{2}(q^2+\zeta^2) & \alpha fb+\beta gq
\end{pmatrix},
\]
which is a lower block triangular matrix. Hence, its eigenvalues are those of the two diagonal $2\times 2$ blocks.

The eigenvalues of $\mathbf{A}_{\mathbf{n}}(\mathbf{u})$ are
\begin{align}\label{gen_eigenvalues}
\tilde{\lambda}_{1,2} = \alpha f b \pm \frac{\alpha f}{2}\, |\mathbf{b}|, \quad \tilde{\lambda}_{3,4}= \sigma \pm \frac{\beta g}{2}\, |\mathbf{q}|, \quad
\tilde{\lambda}_{5,6} = 0,
\end{align}
where
\begin{align*}
|\mathbf{b}| &= \sqrt{b^2+\eta^2}, \qquad
|\mathbf{q}| = \sqrt{q^2+\zeta^2}.
\end{align*}
The corresponding right eigenvectors are given by
\begin{align*}
{\mathbf{\tilde{r}_1}}(\con) &=
\begin{pmatrix}
- f\,\Delta_{-}\\
|\mathbf{b}|\,\Delta_{-}\\
0\\
-\alpha f g\!\left[\delta_{-}(|\mathbf{b}|- b)-\frac{\beta g}{2}\bigl(|\mathbf{b}|q-(bq+\eta\zeta)\bigr)\right]\\
\alpha f\!\left[\frac{\beta g}{2}|\mathbf{q}|^2(|\mathbf{b}|- b)-\delta_{-}\bigl(|\mathbf{b}|q-(bq+\eta\zeta)\bigr)\right]\\
0
\end{pmatrix}, \quad
{\mathbf{\tilde{r}_2}}(\con) =
\begin{pmatrix}
f\,\Delta_{+}\\
|\mathbf{b}|\,\Delta_{+}\\
0\\
-\alpha f g\!\left[\delta_{+}(|\mathbf{b}|+ b)-\frac{\beta g}{2}\bigl(|\mathbf{b}|q+(bq+\eta\zeta)\bigr)\right]\\
\alpha f\!\left[\frac{\beta g}{2}|\mathbf{q}|^2(|\mathbf{b}|+ b)-\delta_{+}\bigl(|\mathbf{b}|q+(bq+\eta\zeta)\bigr)\right]\\
0
\end{pmatrix}, \\
{\mathbf{\tilde{r}_3}}(\con) &=
\begin{pmatrix}
0\\
0\\
0\\
-\dfrac{g}{|\mathbf{q}|}\\
1\\
0
\end{pmatrix}, \qquad
{\mathbf{\tilde{r}_4}}(\con) =
\begin{pmatrix}
0\\
0\\
0\\
\dfrac{g}{|\mathbf{q}|}\\
1\\
0
\end{pmatrix}, \qquad
{\mathbf{\tilde{r}_5}}(\con) =
\begin{pmatrix}
2f\eta\,\Delta_0\\
-4b\eta\,\Delta_0\\
(3b^2-\eta^2)\,\Delta_0\\
\dfrac{\alpha f g}{2}\Bigl(4\alpha b^2\eta f+\beta g(2b\eta q+(3b^2+\eta^2)\zeta)\Bigr)\\
\alpha f\Bigl(\sigma(2b\eta q-(3b^2+\eta^2)\zeta)-\beta b g\eta |\mathbf{q}|^2\Bigr)\\
0
\end{pmatrix}, \\
{\mathbf{\tilde{r}_6}}(\con) &=
\begin{pmatrix}
0\\
0\\
0\\
2\beta^2 g^3\zeta\\
-4\beta g\zeta\,\sigma\\
4\sigma^2-\beta^2 g^2 |\mathbf{q}|^2
\end{pmatrix},
\end{align*}
where
\begin{align*}
\sigma &= \alpha f b + \beta g q,\\
\delta_{\pm} &= \beta g q \pm \frac{\alpha f}{2}\, |\mathbf{b}|, \\
\Delta_{\pm} &= \delta_{\pm}^2 - \frac{\beta^2 g^2}{4}\, |\mathbf{q}|^2, \\
\Delta_0 &= \sigma^2 - \frac{\beta^2 g^2}{4}\, |\mathbf{q}|^2.
\end{align*}
Thus, for every unit vector $\mathbf n\in\mathbb R^2$, the projected one-dimensional system has six real eigenvalues. Moreover, for a fixed direction $\mathbf n$, the repeated eigenvalue $\tilde\lambda_{5,6}=0$ has geometric multiplicity $2$ provided
\[
\Delta_0\neq 0,\qquad f>0,\qquad \mathbf b\neq 0,\qquad 3b^2\neq \eta^2.
\]
On the other hand, a direct computation shows that the four eigenvectors $\tilde{\mathbf r}_i$, $i\in\{1,2,3,4\}$, are linearly independent provided
\[
f>0,\qquad g>0,\qquad \Delta_{\pm}\neq 0,\qquad |\mathbf q|\neq 0,\qquad |\mathbf b|\neq 0.
\]
Hence, for each fixed unit vector $\mathbf n$, the set of eigenvectors $\{\tilde{\mathbf r}_i\}_{i=1}^6$ is linearly independent whenever
\[
f>0,\qquad g>0,\qquad |\mathbf b|\neq 0,\qquad |\mathbf q|\neq 0,\qquad
\Delta_{\pm}\neq 0,\qquad \Delta_0\neq 0,\qquad 3b^2\neq \eta^2.
\]
Therefore, for every fixed direction $\mathbf n$, the projected system is diagonalizable under the above nondegeneracy conditions. Consequently, the two-dimensional system \eqref{eq: Main_system_new_reduced} is hyperbolic in the sense of Definition~\ref{hyper_multi_d} at every state for which the directional matrix $A_{\mathbf n}(\mathbf u)$ is diagonalizable for all unit vectors $\mathbf n\in\mathbb R^2$. 
\subsection{Hyperbolic structure of the system \eqref{eq: Main_system_new_reduced} with zero tangential concentration gradients}
We now consider the special case in which the tangential components of the concentration gradients vanish. Recall that, in the directional system \eqref{eq:directional_system}, these components satisfy
\[
\eta_t=0,\qquad \zeta_t=0.
\]
Hence, if the initial data satisfy
\[
\eta(\cdot,0)=0,\qquad \zeta(\cdot,0)=0,
\]
then
\[
\eta(\cdot,t)=0,\qquad \zeta(\cdot,t)=0
\]
for all $t>0$. Therefore, the manifold
\[
\{\eta=\zeta=0\}
\]
is invariant under the projected dynamics. On this invariant subclass of solutions, the directional system reduces to a $4\times 4$ system of conservation laws involving only the normal variables $(f,b,g,q)$, and the hyperbolic structure becomes considerably simpler.

More precisely, restricting \eqref{eq:directional_system} to the invariant manifold $\{\eta=\zeta=0\}$, we obtain
\begin{equation}\label{1d_system}
\mathbf{u}_t+\mathbf{F}(\mathbf{u})_x=0,
\end{equation}
where
\[
\mathbf{u}=(f,b,g,q)^\top,
\qquad
\mathbf{F}(\mathbf{u})=
\left(
\dfrac{\alpha}{2}f^2b,\;
\dfrac{\alpha}{2}fb^2,\;
\alpha fgb+\dfrac{\beta}{2}g^2q,\;
\dfrac{\beta}{2}gq^2+\alpha fbq
\right)^\top.
\]
The corresponding flux Jacobian is
\[
\mathbf{A}(\mathbf{u})
=
\frac{\partial \mathbf{F}(\mathbf{u})}{\partial \mathbf{u}}
=
\begin{pmatrix}
\alpha fb & \dfrac{\alpha}{2}f^2 & 0 & 0\\[0.3em]
\dfrac{\alpha}{2}b^2 & \alpha fb & 0 & 0\\[0.3em]
\alpha gb & \alpha fg & \alpha fb+\beta gq & \dfrac{\beta}{2}g^2\\[0.3em]
\alpha bq & \alpha fq & \dfrac{\beta}{2}q^2 & \alpha fb+\beta gq
\end{pmatrix}.
\]
In what follows, we consider the system \eqref{eq: Main_system_new_reduced} along the invariant manifold $\{\eta=\zeta=0\}$ and thus consider \eqref{1d_system} as our main system. Moreover, we consider the state space as
\[
\mathcal U
=
\left\{
\mathbf u=(f,b,g,q)^\top\in\mathbb R^4
\;:\;
f>0,\quad b\neq 0,\quad g>0,\quad q>0,\quad
\alpha fb<\beta gq,\quad 2\alpha fb+\beta gq>0
\right\}.
\]
The positivity conditions $f>0$ and $g>0$ select the physically relevant regime, while the additional inequalities
\[
\alpha fb<\beta gq,
\qquad
2\alpha fb+\beta gq>0
\]
ensures the strict hyperbolicity of the system \eqref{1d_system}. We take the assumption $q>0$ to ensure the existence of an explicit Riemann invariant defined in $\mathcal{U}$; see \eqref{selectedinvariants} below.

The eigenvalues $\lambda_k= \lambda_k(\con)$ of $\mathbf{A}(\con)$ can be 
deduced from \eqref{gen_eigenvalues} by setting $\eta=\zeta=0$ and can be ordered for $\con=(f,b,g,q)^{\top}\in \mathcal{U}$ based on the sign of $b$. In particular, for $b>0$ we have the following set of eigenvalues
\begin{align}\label{eigenvalues}
\lambda_1(\con)=\dfrac{\alpha fb}{2}<~\lambda_2(\con)=\dfrac{3\alpha fb}{2}<~\lambda_3(\con)=\alpha fb+\dfrac{\beta gq}{2}<~\lambda_4(\con)=\alpha fb+\dfrac{3\beta gq}{2}.
\end{align}
along with the associated right eigenvectors $\mathbf{r_k}= \mathbf{r_k}(\mathbf{u})$, which reduces to
\begin{equation}\label{righteigenvectors}
\begin{array}{rcl}\displaystyle 
\mathbf{r_1}(\con)= \left(-\frac{f}{b},1,0,0\right)^{\top}, && \!\!\!\!\!\! \displaystyle \mathbf{r_2}(\con)=\left(\frac{\alpha fb-3\beta gq}{4\alpha bq},\frac{\alpha fb-3\beta gq}{4\alpha fq},\frac{g}{q},1\right)^{\top},\\[3ex]\displaystyle 
\mathbf{r_3}(\con)=\left(0,0,-\frac{g}{q},1\right)^{\top}, & & \!\!\!\!\!\! \displaystyle \mathbf{r_4}(\con)=\left(0,0,\frac{g}{q},1\right)^{\top}.
\end{array}
\end{equation}
Similarly, for $b<0$, we have the following ordering of eigenvalues for $\mathbf{u}\in \mathcal{U}$
\begin{align}\label{eigenvalues}
\lambda_1(\con)=\dfrac{3\alpha fb}{2}<~\lambda_2(\con)=\dfrac{\alpha fb}{2}<~\lambda_3(\con)=\alpha fb+\dfrac{\beta gq}{2}<~\lambda_4(\con)=\alpha fb+\dfrac{3\beta gq}{2}.
\end{align}
The eigenvectors are linearly independent in $\mathcal U$.  This implies the hyperbolicity of the reduced system \eqref{1d_system} in $\mathcal U$. 
\subsection{Riemann Invariants of the system \eqref{1d_system}}
Using the Definition \ref{Defn_RI}, one can obtain the complete set of Riemann invariants for the system \eqref{1d_system} associated with each
characteristic field. The explicit Riemann invariants are tabulated in Table~\ref{tab:RI_fields}.
\begin{table}[!ht]
\centering
\begin{tabular}{c|c}
\hline
Characteristic Speed & Independent Riemann Invariants \\
\hline

$\lambda_1=\dfrac{\alpha fb}{2}$ 
& $fb,\; g,\; q$ \\[1.5ex]

$\lambda_2=\dfrac{3\alpha fb}{2}$ 
& $\dfrac{b}{f},\; \dfrac{q}{g},\; \dfrac{\alpha fb+\beta gq}{\alpha(gq)^{1/4}}$ \\[1.5ex]

$\lambda_3=\alpha fb+\dfrac{\beta gq}{2}$ 
& $f,\; b,\; gq$ \\[1.5ex]

$\lambda_4=\alpha fb+\dfrac{3\beta gq}{2}$ 
& $f,\; b,\; \dfrac{g}{q}$ \\

\hline
\end{tabular}
\caption{\vspace{0.3 cm} All independent Riemann invariants associated with each characteristic field.}
\label{tab:RI_fields}
\end{table}

Selecting one invariant from each characteristic family in
Table~\ref{tab:RI_fields}, we define the transformation
\[
\con=(f,b,g,q)^T
\;\longmapsto\;
\mathbf{V}=(v_1,v_2,v_3,v_4)^T,
\]
with
\begin{equation}\label{selectedinvariants}
v_1=\frac{b}{f},
\qquad
v_2=fb,
\qquad
v_3=\frac{q}{g},
\qquad
v_4=\frac{\alpha fb+\beta gq}{\alpha(gq)^{1/4}}.
\end{equation}

This choice yields a complete coordinate system of Riemann invariants,
and therefore fully diagonalizes the system in the $\mathbf{V}$–variables.
The transformed system reads
\begin{equation}\label{Rtransformed}
\begin{array}{rclrcl}
(v_1)_t+\dfrac{\alpha v_2}{2}(v_1)_x&=&0,\\[1ex]
(v_2)_t+\dfrac{3\alpha v_2}{2}(v_2)_x&=&0,\\[1ex]
(v_3)_t+\left(\alpha v_2+\dfrac{\beta \psi(v_1, v_4)}{2}\right)(v_3)_x&=&0,\\[1ex]
(v_4)_t+\left(\alpha v_2+\dfrac{3\beta \psi(v_1, v_4)}{2}\right)(v_4)_x&=&0,
\end{array}
\end{equation}
where $\psi(v_1, v_4)=gq$, which can be obtained by solving the nonlinear quartic equation
\[
\alpha v_2+\beta x-\alpha x^{1/4} v_4=0.
\]
Though this transformation fully diagonalizes the system, the inversion $\mathbf{V} \mapsto \con$ is algebraically complex and requires a quartic equation to be solved, which turns out to be computationally expensive while being less accurate due to the iterative solvers required for the inversion process; see test problem 1 in Section \ref{Ex1}. To get a more computationally efficient transformation, we select different Riemann invariants from Table~\ref{tab:RI_fields} that still yield a sparse structure in the transformed eigenvector matrices but allow for a simpler inversion. In particular, we select the following Riemann invariants:
\begin{equation}\label{selectedinvariants2}
v_1 = \frac{b}{f},
\quad
v_2 = fb,
\quad
v_3 = \frac{q}{g},
\quad
v_4 = gq.
\end{equation}

The normalized transformed right and left eigenvector matrices then take the following sparse form:

\begin{align}\label{transformed_eigenvectors}
    R =
    \begin{pmatrix}
    1 & 0 & 0 & 0 \\
    0 & 1 & 0 & 0 \\
    0 & 0 & 1 & 0 \\
    0 & \dfrac{4 \alpha g q}{\alpha b f - 3\beta g q} & 0 & 1
    \end{pmatrix},
    \quad 
    L =
    \begin{pmatrix}
    1 & 0 & 0 & 0 \\
    0 & 1 & 0 & 0 \\
    0 & 0 & 1 & 0 \\
    0 & -\dfrac{4 \alpha g q}{\alpha b f - 3\beta g q} & 0 & 1
    \end{pmatrix}.
\end{align}
It is worth noting that, for the system \eqref{eq: Main_system_old}, the entropy--entropy flux pairs derived in \cite{barthwal2025hyperbolic} were constructed using the fact that the Riemann invariants form a coordinate system. In particular, the diagonal structure of the system was exploited to identify a broad class of entropy functions; see \cite{barthwal2025hyperbolic} for more details. A similar approach can also be employed in the present setting to derive explicit entropy/entropy-flux pairs, including a strictly convex entropy. We provide one such strictly convex entropy of the hyperbolic system \eqref{1d_system} without proof. The entropy of the hyperbolic system \eqref{1d_system} takes the form 
\[
E(f, b, g, q)=\dfrac{1}{|fb|}+\dfrac{1}{\alpha fb+\beta gq}
\]
Consequently, the local well-posedness of classical solutions for the Cauchy problem for the system \eqref{1d_system} can be established using the classical theory of hyperbolic balance laws. However, in order to keep the focus of the present article on the computational aspects, we do not pursue the entropy structure and associated nonlinear waves of \eqref{1d_system} in detail and instead refer the interested readers to \cite{barthwal2025hyperbolic} for a more intensive analysis.
\section{Numerical Methodology}\label{sec: Numerical_scheme}
We are interested in developing a high-order finite difference WENO scheme for the hyperbolic system \eqref{1d_system}. We adopt the alternative WENO (AWENO) formulation proposed in \cite{wu2025finite}, which uses a Riemann invariant-based local characteristic decomposition (RI-LCD) to achieve high-order accuracy and non-oscillatory properties with a more computationally efficient implementation compared to the classical WENO schemes.

We start by considering the one-dimensional system \eqref{1d_system}. We discretize the computational domain with a uniform grid $I_j = [x_{j-1/2}, x_{j+1/2}]$ with grid size $\Delta x$. Then a semi-discrete finite difference scheme for \eqref{1d_system} can be written as
\begin{equation}\label{semi_discrete_scheme}
    \frac{d}{dt} \con_j + \frac{1}{\Delta x} \left( \hat{\mathbf{F}}_{j+1/2} - \hat{\mathbf{F}}_{j-1/2} \right) = \mathbf{0},
\end{equation}
where $\con_j$ approximates $\con(t, x_j)$ and $\hat{\mathbf{F}}_{j+1/2}$ is the numerical flux at the cell interface $x_{j+1/2}$ and $k$-th order accurate in the sense that
\begin{align*}
    \frac{1}{\Delta x} \left( \hat{\mathbf{F}}_{j+1/2} - \hat{\mathbf{F}}_{j-1/2} \right) = \mathbf{F}(\con(x))_x|_{x_j} + \mathcal{O}(\Delta x^k).
\end{align*}

\subsection{Alternative WENO Formulations}
Unlike classical finite difference WENO schemes, which reconstruct split fluxes, the AWENO scheme operates on the point values of conserved variables. The numerical flux $\hat{\mathbf{F}}_{j+1/2}$ is computed via the following steps \cite{jiang2013alternative,balsara2025efficient}:
\begin{enumerate}
    \item \textbf{Preparation:} Obtain ghost point values via boundary conditions. Compute the flux function values $\mathbf{F}_j = \mathbf{F}(\con_j)$ for all grid points.

    \item \textbf{WENO Interpolation:} Apply the WENO interpolation to the conserved variables $\con$ component-wise to obtain high-order reconstructed values $\con^{\pm}_{j+1/2}$ at the interface to obtain
    \begin{align}
        \con^-_{j+1/2} &= \text{WENO}(\con_{j-r+1}, \dots, \con_{j+r-1}),\\
        \con^+_{j+1/2} &= \text{WENO}(\con_{j+r}, \dots, \con_{j-r+2}).
    \end{align}

    \item \textbf{Flux Assembly:} The final numerical flux is the sum of a low-order flux and a high-order correction:
    \begin{equation}
        \hat{\mathbf{F}}_{j+1/2} = \hat{\mathbf{F}}^{\text{low}}_{j+1/2} + \hat{\mathbf{F}}^{\text{cor}}_{j+1/2}.
    \end{equation}
    Here, $\hat{\mathbf{F}}^{\text{low}}_{j+1/2} = \mathcal{F}(\con^-_{j+1/2}, \con^+_{j+1/2})$ is computed using an approximate Riemann solver (e.g., HLL flux). The high-order correction term $\hat{\mathbf{F}}^{\text{cor}}_{j+1/2}$ is computed directly from the flux values $\mathbf{F}_j$.
    \end{enumerate}

The high-order flux correction terms $\hat{\mathbf{F}}^{\text{cor}}_{j+1/2}$ ensure the scheme achieves the desired order of accuracy $k = 2r-1$. The coefficients for various orders can be found in \cite{balsara2025efficient}. In this work, we will use the following correction terms for the 5th order scheme, which corresponds to $r=3$, and given by
    \begin{equation}
        \hat{\mathbf{F}}^{\text{cor},3}_{j+1/2} = \frac{19}{3840}(\mathbf{F}_{j-2} + \mathbf{F}_{j+3}) - \frac{137}{3840}(\mathbf{F}_{j-1} + \mathbf{F}_{j+2}) + \frac{59}{1920}(\mathbf{F}_{j} + \mathbf{F}_{j+1}).
    \end{equation}

If we apply WENO interpolation to the conservative variable $\con$, the resulting scheme is called the \emph{Componentwise WENO} (CP-WENO) scheme. This is the simplest approach but can lead to spurious oscillations near discontinuities which we will demonstrate in the numerical results section. We will now present two more approaches that are designed to mitigate these oscillations by reconstructing in the characteristic space.
\subsection{Local Characteristic Decomposition}\label{LCD}
In the LCD (local characteristic decomposition) approach, the characteristic decomposition is performed locally at each interface $x_{j+\frac12}$. We first compute the Jacobian matrix of the flux function and evaluate it at the interface using locally averaged states. In our implementation, we use simple arithmetic averaging. In particular,
$$
\con_{j+\frac12} = \frac{\con_j + \con_{j+1}}{2}.
$$

Using this interface state, the Jacobian is diagonalized as
$$
\mathbf{A}_{j+\frac12} = \mathbf{A}(\con_{j+\frac12})= \mathbf{R}_{j+\frac12}\,\Lambda_{j+\frac12}\,\mathbf{L}_{j+\frac12},
$$
where $\mathbf{R}_{j+\frac12}$ and $\mathbf{L}_{j+\frac12}$ denote the right and left eigenvector matrices, respectively, and $\Lambda_{j+\frac12}$ is the diagonal matrix of eigenvalues.

The conservative variables in the reconstruction stencil are then projected onto the local characteristic space:

\begin{align}\label{u_to_w}
    \mathbf{w}_k = \mathbf{L}_{j+\frac12} \con_k, 
    \qquad k = j-r+1,\ldots,j+r.
\end{align}

The WENO interpolation is then carried out component-wise on the characteristic variables as 
\begin{align}\label{w_to_u}
        \mathbf{w}^-_{j+1/2} &= \text{WENO}(\mathbf{w}_{j-r+1}, \dots, \mathbf{w}_{j+r-1}),\\
        \mathbf{w}^+_{j+1/2} &= \text{WENO}(\mathbf{w}_{j+r}, \dots, \mathbf{w}_{j-r+2}).
\end{align}
These are finally transformed back to physical space via
$$
\con_{j+\frac12}^{\pm} = \mathbf{R}_{j+\frac12} \mathbf{w}_{j+\frac12}^{\pm}.
$$

\subsection{Riemann Invariant based Local Characteristic Decomposition} \label{RI-LCD}
The LCD approach removes spurious oscillations from the solution but is computationally expensive due to the need for several matrix-vector products in \eqref{u_to_w} and \eqref{w_to_u}. To mitigate this cost, we can use the Riemann invariant-based LCD (RI-LCD) approach. Instead of projecting onto the local characteristic space, we perform the WENO interpolation directly on the Riemann invariants defined in \eqref{selectedinvariants2}. This allows us to achieve a sparse structure in the transformed eigenvector matrices, which significantly reduces the computational cost while still effectively suppressing oscillations. The following steps summarize the RI-LCD procedure:
\begin{enumerate}
    \item \textbf{Variable Transformation:} Transform the conservative variables $\con_j$ to the Riemann invariant variables $\mathbf{V}_j$ using the transformation defined in \eqref{selectedinvariants2}.

    \item \textbf{WENO Interpolation:} Apply the WENO interpolation to the Riemann invariant variables $\mathbf{V}$ using the LCD approach presented in Section \ref{LCD} with the transformed eigenvector matrices defined in \eqref{transformed_eigenvectors} to obtain $\mathbf{V}^\pm_{j+1/2}$ .

    \item \textbf{Inverse Transformation:} Transform the reconstructed Riemann invariant values $\mathbf{V}^\pm_{j+1/2}$ back to conservative variables using the inverse of the transformation defined in \eqref{selectedinvariants2}.
\end{enumerate}
\subsection{Time Integration}
For time integration of the semi-discrete scheme, we use the following third-order strong stability preserving Runge-Kutta (SSP-RK3) method \cite{gottlieb2001strong}:
\begin{subequations}\label{SSP-RK3}
\begin{align}
\con^{(1)} &= \con^n + \Delta t \, L(\con^n), \\
\con^{(2)} &= \frac{3}{4} \con^n + \frac{1}{4} \con^{(1)} + \frac{1}{4} \Delta t \, L(\con^{(1)}), \\
\con^{n+1} &= \frac{1}{3} \con^n + \frac{2}{3} \con^{(2)} + \frac{2}{3} \Delta t \, L(\con^{(2)}).
\end{align}
\end{subequations}
where $L(\con)$ represents the spatial discretization operator defined by the right-hand side
\begin{align*}
L(\con_j) = -\frac{1}{\Delta x} \left( \hat{\mathbf{F}}_{j+1/2} - \hat{\mathbf{F}}_{j-1/2} \right).
\end{align*}
Time step $\Delta t$ is chosen according to the CFL condition
\begin{align*}
\Delta t \leq \text{CFL} \cdot \frac{\Delta x}{\max_j |\lambda_{\max}(\con_j)|}.
\end{align*}
In one exception, for the accuracy verification test problem, we use $ \Delta t = \Delta x^{\frac{5}{3}}$
to ensure that the temporal error does not dominate the spatial error, allowing us to accurately measure the convergence rates of the schemes. 
\subsection{Two dimensional Extension}
The proposed schemes can be extended to two dimensions using a dimension-by-dimension approach. Consider the two-dimensional system
\begin{align*}
\con_t + \mathbf{F}(\con)_x + \mathbf{G}(\con)_y = 0.
\end{align*}
We can discretize the spatial domain with uniform grid cells $I_{i,j} = [x_{i-1/2}, x_{i+1/2}] \times [y_{j-1/2}, y_{j+1/2}]$ with grid sizes $\Delta x$ and $\Delta y$. The semi-discrete scheme can be written as
\begin{align*}
\frac{d}{dt} \con_{i,j} + \frac{1}{\Delta x} \left( \hat{\mathbf{F}}_{i+1/2,j} - \hat{\mathbf{F}}_{i-1/2,j} \right) + \frac{1}{\Delta y} \left( \hat{\mathbf{G}}_{i,j+1/2} - \hat{\mathbf{G}}_{i,j-1/2} \right) = \mathbf{0},
\end{align*}
where $\hat{\mathbf{F}}_{i+1/2,j}$ and $\hat{\mathbf{G}}_{i,j+1/2}$ are the numerical fluxes in the $x$ and $y$ directions, respectively. The numerical fluxes can be computed using the same AWENO procedure as described in the one-dimensional case, applied separately in each spatial direction. The time integration can be performed using the same SSP-RK3 method as in the one-dimensional case \eqref{SSP-RK3}, with
\begin{align*}
L(\con_{i,j}) = -\frac{1}{\Delta x} \left( \hat{\mathbf{F}}_{i+1/2,j} - \hat{\mathbf{F}}_{i-1/2,j} \right) - \frac{1}{\Delta y} \left( \hat{\mathbf{G}}_{i,j+1/2} - \hat{\mathbf{G}}_{i,j-1/2} \right),
\end{align*}
and the time step $\Delta t$ chosen as

\begin{align*}
    \Delta t \leq \text{CFL} \cdot \frac{1}{\frac{\lambda_{\max}^x}{\Delta x} + \frac{\lambda_{\max}^y}{\Delta y}},
\end{align*}

where $\lambda_{\max}^x$ and $\lambda_{\max}^y$ are the maximum wave speeds in the $x$ and $y$ directions, respectively.
\section{Numerical Results}
In this section, we present a set of numerical experiments to evaluate the performance of the proposed schemes. We compare the Componentwise WENO (CP-WENO), Local Characteristic Decomposition WENO (LCD-WENO), Riemann Invariant-based Local Characteristic Decomposition WENO (RI-WENO), and Complete set of Riemann invariant-based Local Characteristic Decomposition WENO (CRI-WENO) schemes for several test problems involving both smooth and discontinuous solutions. Whenever possible, analytical solutions are provided for reference, and the solution by the Godunov scheme is included as a baseline for comparison.

We have used the fifth-order WENO interpolation for all the schemes in the numerical experiments. The CFL is set to $0.4$ for all simulations.
\subsection{One-dimensional test cases}
In this section, we take some one-dimensional test cases from \cite{barthwal2025hyperbolic} and \cite{barthwal2026generalized} and validate the numerical results using the scheme developed in Section \ref{sec: Numerical_scheme}. To investigate the effect of different Marangoni numbers on the evolution of the film heights, we also provide a new test case.
\testproblem{Accuracy and efficiency test}\label{Ex1}
To verify the formal order of accuracy of the scheme, we consider a smooth
manufactured solution with periodic boundary conditions in the $x$-direction from \cite{barthwal2026generalized}. The exact solution is given by
$$
f(x,t) = 2 + \sin(x+t), \quad
b(x,t) = -\frac{2}{f(x,t)}, \quad
g(x,t) = 2, \quad
q(x,t) = 1.
$$
The computational domain is $x \in [-\pi, \pi]$. The final simulation time is taken as $ t = 1$. We have presented the $L^1$, $L^2$, and $L^\infty$ errors and convergence rates for all three schemes in Table~\ref{tab:accuracy}. The results confirm that all schemes achieve the expected fifth-order accuracy for smooth solutions. In addition, the RI-WENO scheme demonstrates slightly smaller errors compared to the CP-WENO, LCD-WENO, and CRI-WENO schemes, which can be attributed to the more efficient handling of the characteristic structure of the system.

\begin{figure}[ht]
    \centering
    \begin{tabular}{cc}
        \includegraphics[width=0.48\textwidth]{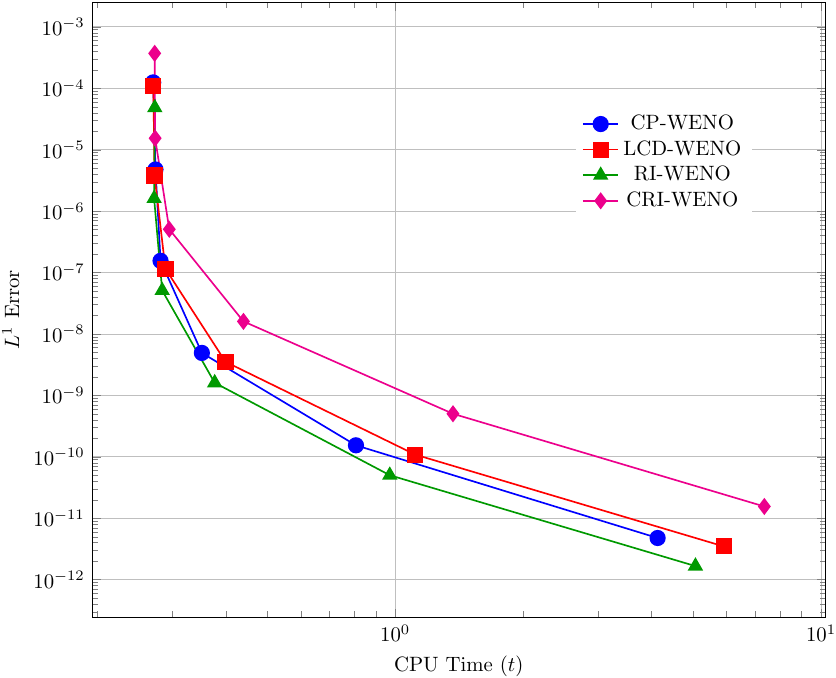} &
        \includegraphics[width=0.48\textwidth]{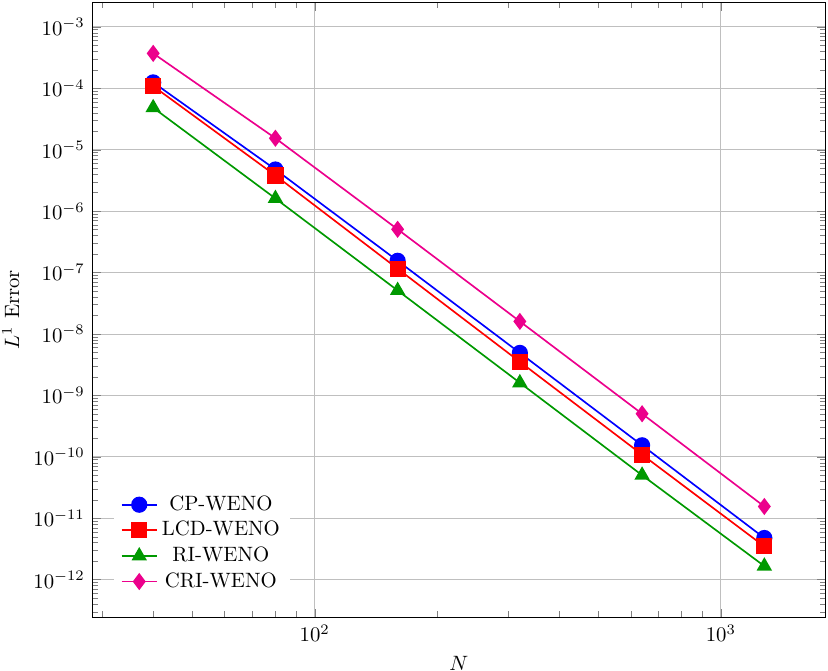} \\
        (a) CPU time versus $L^1$-error & (b)  $N$ versus $L^1$-error.
    \end{tabular}
    \caption{Comparison plots for Test Problem \ref{Ex1} of RI-WENO, CP-WENO, LCD-WENO and CRI-WENO schemes. }
    \label{fig:l1plots}
\end{figure}

In Figure~\ref{fig:l1plots}, we compare the performance of the considered schemes in terms of CPU time required to achieve the $L^1$ error on a log--log scale, along with the corresponding convergence plots on a log--log scale. From Figure~\ref{fig:l1plots} (a), it can be observed that the RI-WENO scheme requires significantly less CPU time to attain the same $L^1$ error compared to the other schemes. Furthermore, Figure~\ref{fig:l1plots} (b) shows that the RI-WENO scheme provides better accuracy than the other considered schemes. In particular, the scheme based on the full set of Riemann invariants, namely CRI-WENO, is less accurate and computationally more expensive because of the iterative solver required for the inversion process. This observation indicates that, even when a system possesses a full set of Riemann invariants, the computational cost may still be high if the transformation back to the original variables is algebraically complex. It also suggests that, although Riemann-invariant-based methods can improve accuracy, careful selection of a suitable set of Riemann invariants is essential, rather than relying directly on the set of Riemann invariants forming a coordinate system. 
\begin{table}[H]
\centering
\caption{Convergence study for CP-WENO, LCD-WENO, RI-WENO, and CRI-WENO schemes.\vspace{0.3cm}}
\label{tab:accuracy}
\begin{tabular}{c c c c c c c}
\toprule
$N$ & $L^1$ error & $L^1$ rate & $L^2$ error & $L^2$ rate & $L^\infty$ error & $L^\infty$ rate \\
\midrule

\multicolumn{7}{c}{\textbf{CP-WENO}} \\
\midrule
20  & 2.378658e-03 & --    & 2.817288e-03 & --    & 5.589221e-03 & --    \\
40  & 1.247322e-04 & 4.253 & 1.718981e-04 & 4.035 & 4.212099e-04 & 3.730 \\
80  & 4.763562e-06 & 4.711 & 7.066511e-06 & 4.604 & 2.268759e-05 & 4.215 \\
160 & 1.553167e-07 & 4.939 & 2.369187e-07 & 4.899 & 8.659043e-07 & 4.712 \\
320 & 4.914033e-09 & 4.982 & 7.499654e-09 & 4.981 & 2.834420e-08 & 4.933 \\
640 & 1.539123e-10 & 4.997 & 2.345452e-10 & 4.999 & 8.909335e-10 & 4.992 \\

\addlinespace
\multicolumn{7}{c}{\textbf{LCD-WENO}} \\
\midrule
20  & 3.086493e-03 & --    & 3.770919e-03 & --    & 8.741611e-03 & --    \\
40  & 1.084724e-04 & 4.831 & 1.441619e-04 & 4.709 & 3.213477e-04 & 4.766 \\
80  & 3.814721e-06 & 4.830 & 5.608207e-06 & 4.684 & 1.652505e-05 & 4.281 \\
160 & 1.155005e-07 & 5.046 & 1.791501e-07 & 4.968 & 6.219107e-07 & 4.732 \\
320 & 3.521138e-09 & 5.036 & 5.557993e-09 & 5.010 & 2.049343e-08 & 4.923 \\
640 & 1.085128e-10 & 5.020 & 1.722960e-10 & 5.012 & 6.469862e-10 & 4.985 \\

\addlinespace
\multicolumn{7}{c}{\textbf{RI-WENO}} \\
\midrule
20  & 1.453345e-03 & --    & 1.731971e-03 & --    & 3.585057e-03 & --    \\
40  & 4.830865e-05 & 4.911 & 5.984330e-05 & 4.855 & 1.453204e-04 & 4.625 \\
80  & 1.608455e-06 & 4.909 & 1.919781e-06 & 4.962 & 4.628580e-06 & 4.973 \\
160 & 5.108325e-08 & 4.977 & 6.010381e-08 & 4.997 & 1.396106e-07 & 5.051 \\
320 & 1.603981e-09 & 4.993 & 1.868474e-09 & 5.008 & 4.146463e-09 & 5.073 \\
640 & 5.016161e-11 & 4.999 & 5.813561e-11 & 5.006 & 1.234901e-10 & 5.069 \\

\addlinespace
\multicolumn{7}{c}{\textbf{CRI-WENO}} \\
\midrule
20  & 5.322221e-03 & --    & 6.957962e-03 & --    & 1.387233e-02 & --    \\
40  & 3.705585e-04 & 3.844 & 5.581100e-04 & 3.640 & 1.286650e-03 & 3.431 \\
80  & 1.537824e-05 & 4.591 & 2.491730e-05 & 4.485 & 7.390992e-05 & 4.122 \\
160 & 5.063395e-07 & 4.925 & 8.549406e-07 & 4.865 & 3.125959e-06 & 4.563 \\
320 & 1.602438e-08 & 4.982 & 2.719388e-08 & 4.974 & 1.068161e-07 & 4.871 \\
640 & 5.021174e-10 & 4.996 & 8.508649e-10 & 4.998 & 3.413641e-09 & 4.968 \\

\bottomrule
\end{tabular}
\end{table}
\testproblem{Large rarefaction test case\label{RP9}}
We consider a one-dimensional Riemann problem from \cite{barthwal2026generalized}. The domain is $x \in [-20,5]$ with outflow boundary conditions. The initial condition consists of two constant states separated by a discontinuity at $x=0$, given by

$$
(f,b,g,q)(x,0)=
\begin{cases}
(2.0,\,-2.0,\,16.0,\,16.0/7.0), & x \le 0, \\
(1.0,\,-1.0,\,4.0,\,4.0/7.0), & x > 0.
\end{cases}
$$

\begin{figure}[htbp]
\centering
\begin{subfigure}{0.49\textwidth}
\includegraphics[width=\linewidth]{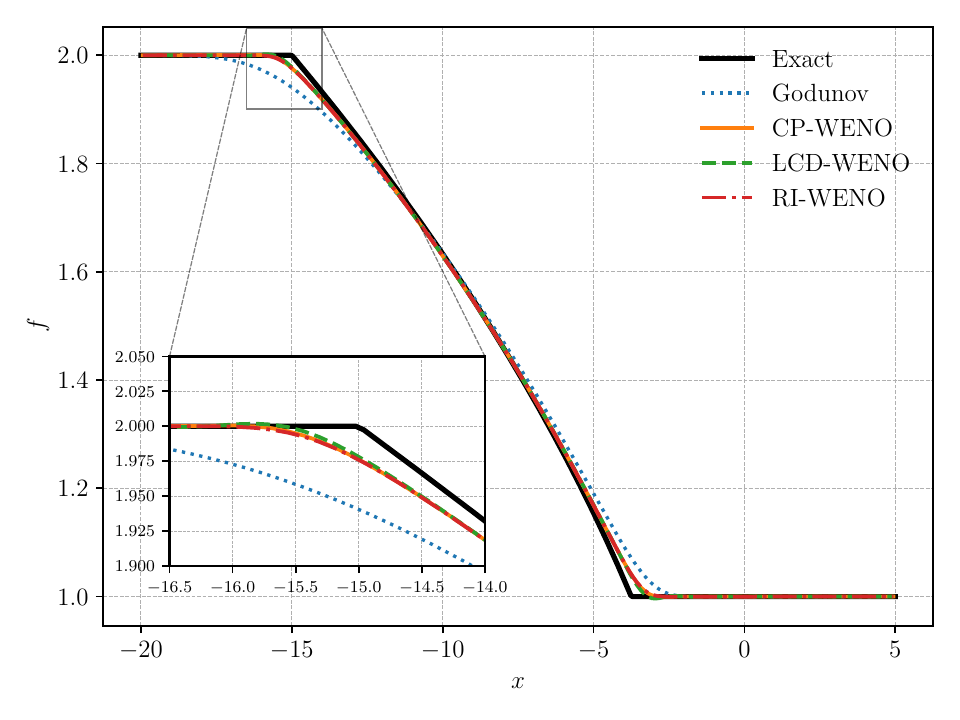}
\caption{Film height ($f$)}
\end{subfigure}\hspace{0pt}
\begin{subfigure}{0.49\textwidth}
\includegraphics[width=\linewidth]{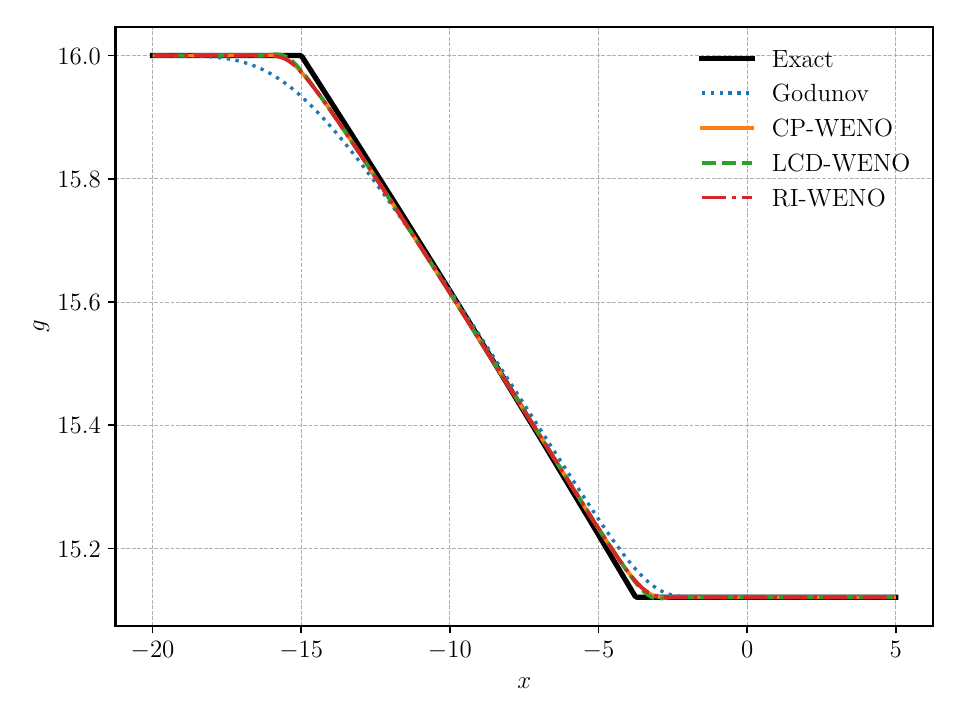}
\caption{Film height ($g$)}
\end{subfigure}

\medskip

\begin{subfigure}{0.49\textwidth}
\includegraphics[width=\linewidth]{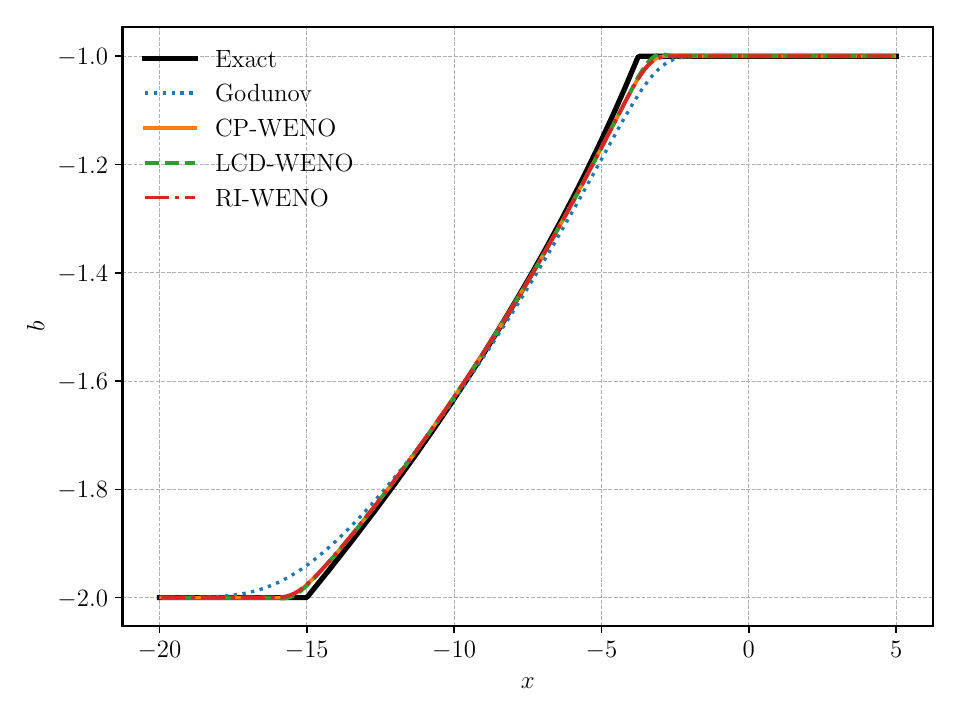}
\caption{Concentration gradient ($b$)}
\end{subfigure}\hspace{0pt}
\begin{subfigure}{0.49\textwidth}
\includegraphics[width=\linewidth]{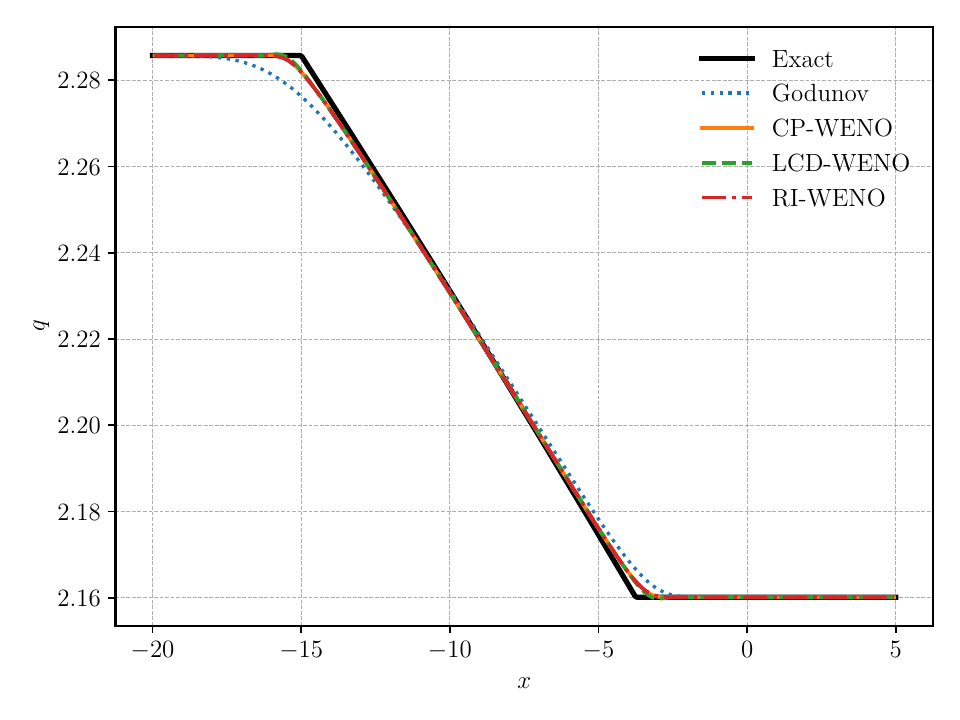}
\caption{Concentration gradient ($q$)}
\end{subfigure}
\caption{Comparison of the CP-WENO, LCD-WENO, and RI-WENO schemes with the exact and Godunov solutions for Test Problem~\ref{RP9} on a mesh with $N_x=200$ at the final time $t=2.5$.}
\label{fig:RP9_all}
\end{figure}

The problem is computed on a uniform mesh with $N_x=200$ cells, and the solution is evolved up to the final time $t=2.5$. The numerical solutions obtained using different schemes are compared in Figure~\ref{fig:RP9_all}. The results show that the CP-WENO, LCD-WENO, and RI-WENO schemes successfully capture the large rarefaction wave and yield very similar solution profiles. We also provide a CPU time comparison for different Riemann data in Table \ref{tab:cpu_time}. It can be observed that the RI-WENO scheme outperforms its counterparts in the CPU time here as well.

\testproblem{Shock tube type test case\label{RP2}}

\begin{figure}[htbp]
\centering
\begin{subfigure}{0.48\textwidth}
\includegraphics[width=\linewidth]{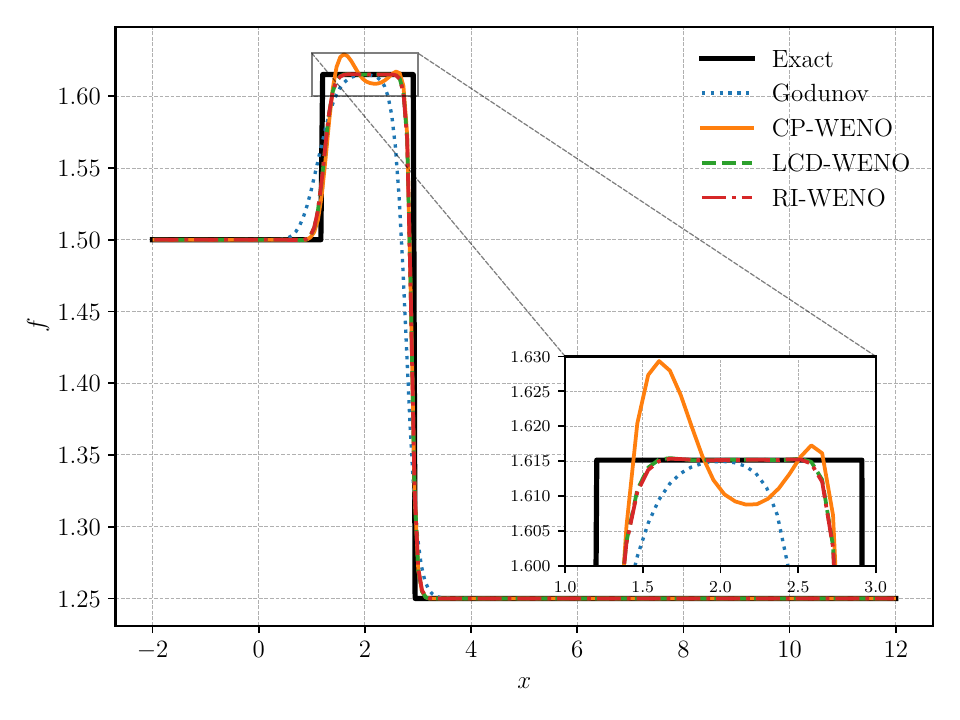}
\caption{Film height ($f$)}
\end{subfigure}
\hfill
\begin{subfigure}{0.48\textwidth}
\includegraphics[width=\linewidth]{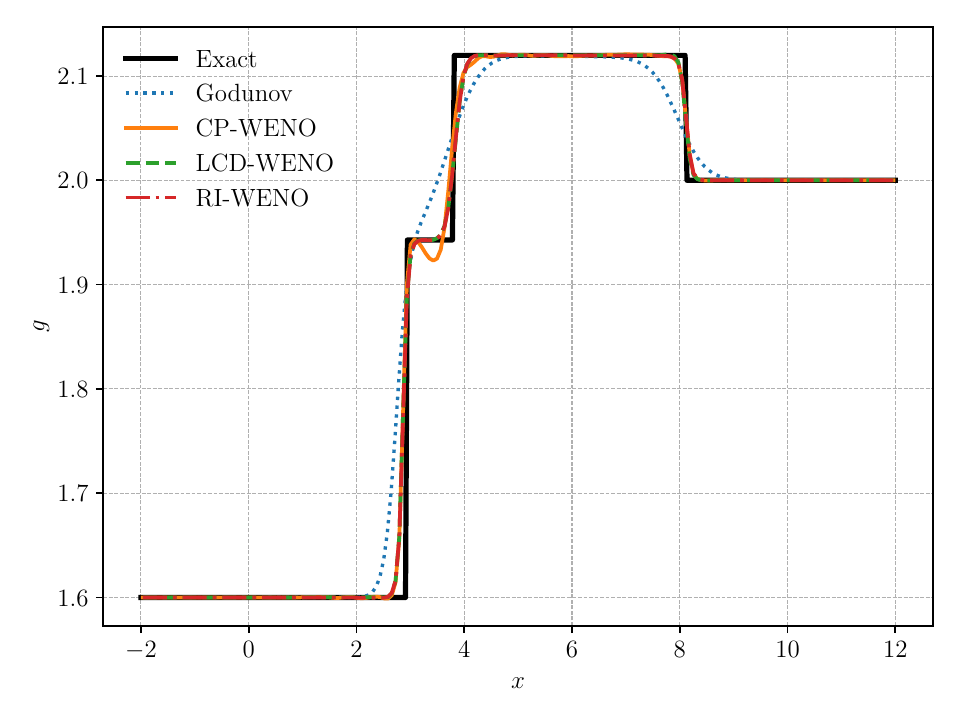}
\caption{Film height ($g$)}
\end{subfigure}

\vspace{0.3cm}

\begin{subfigure}{0.48\textwidth}
\includegraphics[width=\linewidth]{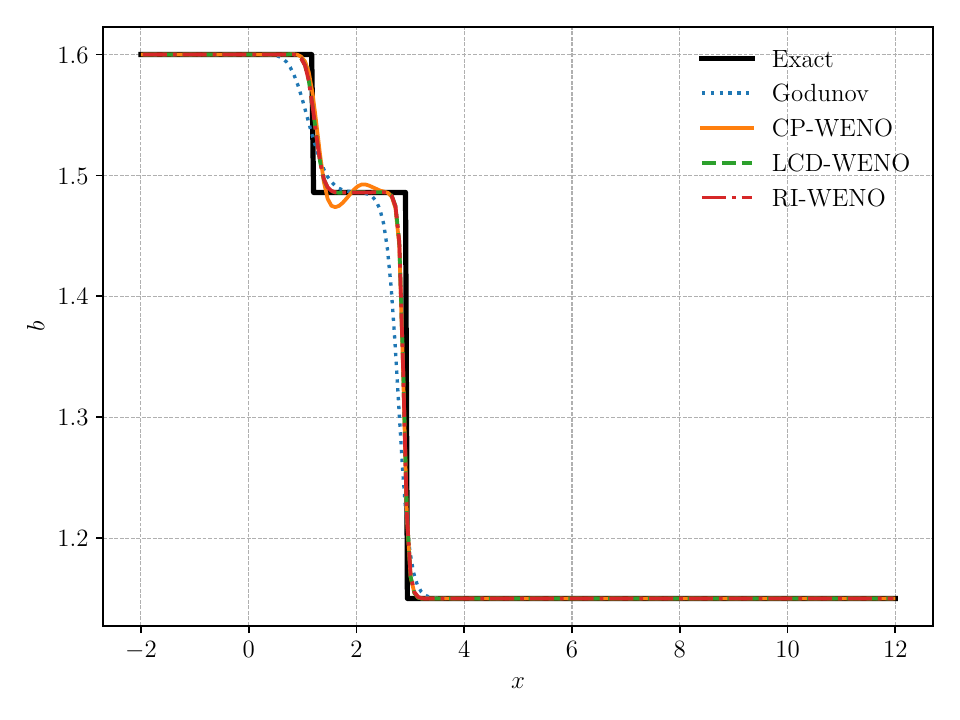}
\caption{Concentration gradient ($b$)}
\end{subfigure}
\hfill
\begin{subfigure}{0.48\textwidth}
\includegraphics[width=\linewidth]{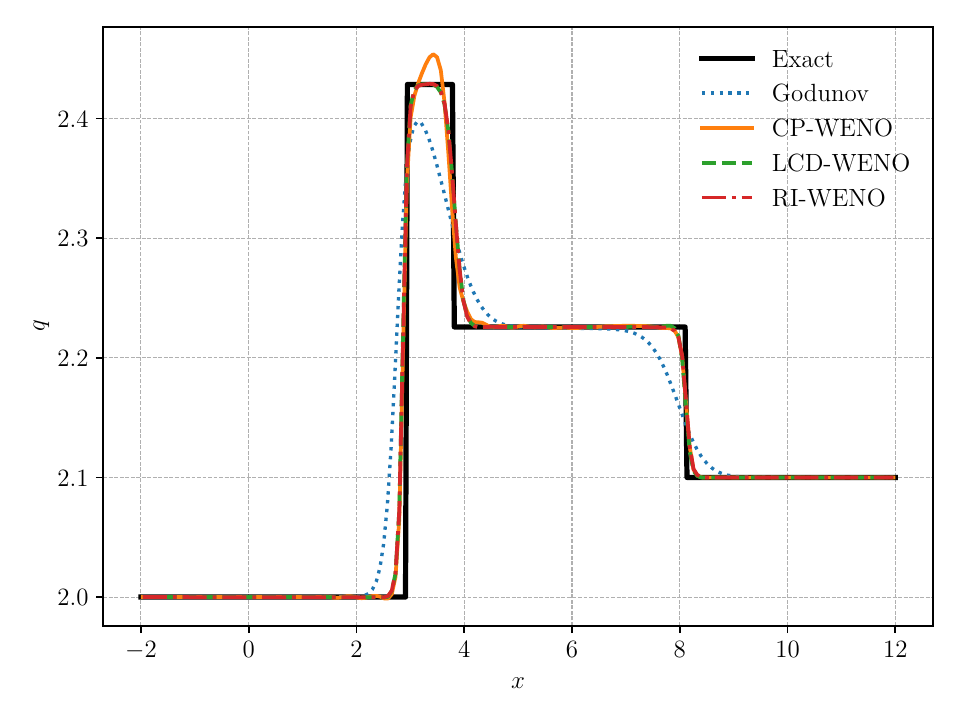}
\caption{Concentration gradient ($q$)}
\end{subfigure}

\caption{Comparison of CP-WENO, LCD-WENO, and RI-WENO with the exact and Godunov solutions for Test Problem~\ref{RP2} at $N_x=200$ and $t=1.0$.}
\label{fig:RP2_all}
\end{figure}

We choose the second Riemann problem from \cite{barthwal2025hyperbolic} to test the performance of the schemes in capturing shock waves and contact discontinuities. We fill the initial condition within the domain $x \in [-2,12]$ with two constant states separated at $x=0$, given by
$$
(f,b,g,q)(x,0)=
\begin{cases}
(1.5,\,1.6,\,1.6,\,2.00), & x \le 0, \\
(1.25,\,1.15,\,2.0,\,2.10), & x > 0.
\end{cases}
$$
We impose outflow boundary conditions on both ends of the domain. The final time is taken as $t=1.0$. We compute the solution on a uniform grid with $N_x=200$ cells. The numerical solutions obtained using CP-WENO, LCD-WENO, and RI-WENO are compared against the exact solution and the Godunov scheme in Figure~\ref{fig:RP2_all}. The CP-WENO scheme exhibits spurious oscillations near the discontinuity, while both LCD-WENO and RI-WENO effectively suppress these oscillations. The RI-WENO scheme achieves accuracy comparable to LCD-WENO while being computationally more efficient, as demonstrated by the CPU time comparison in Table~\ref{tab:cpu_time}.

\begin{figure}[htbp]
\centering
\begin{subfigure}{0.48\textwidth}
\includegraphics[width=\linewidth]{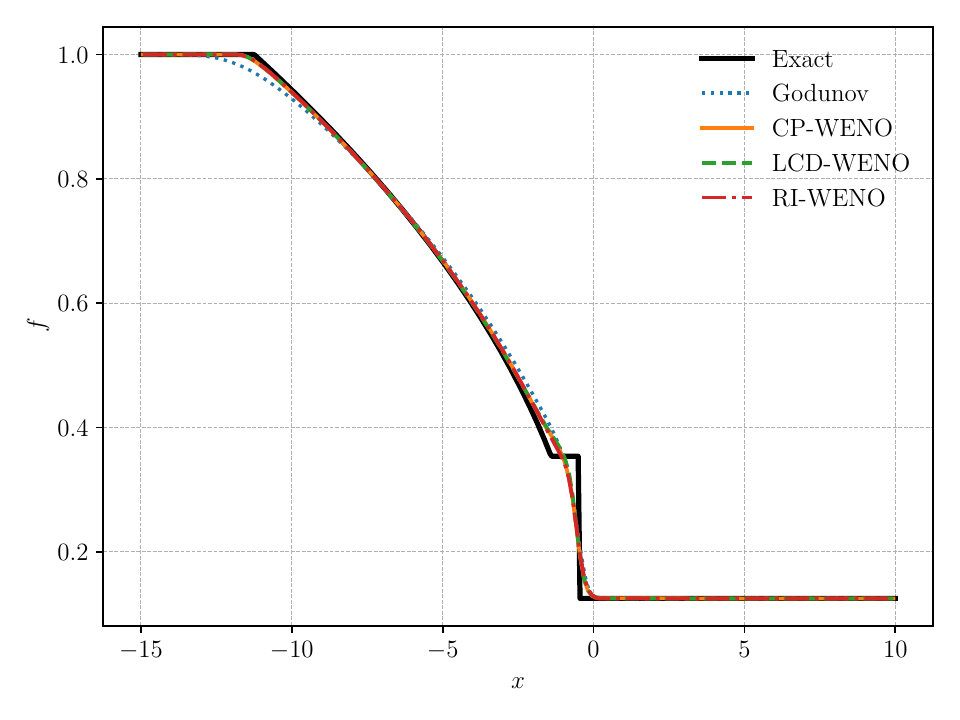}
\caption{Film height ($f$)}
\end{subfigure}
\hfill
\begin{subfigure}{0.48\textwidth}
\includegraphics[width=\linewidth]{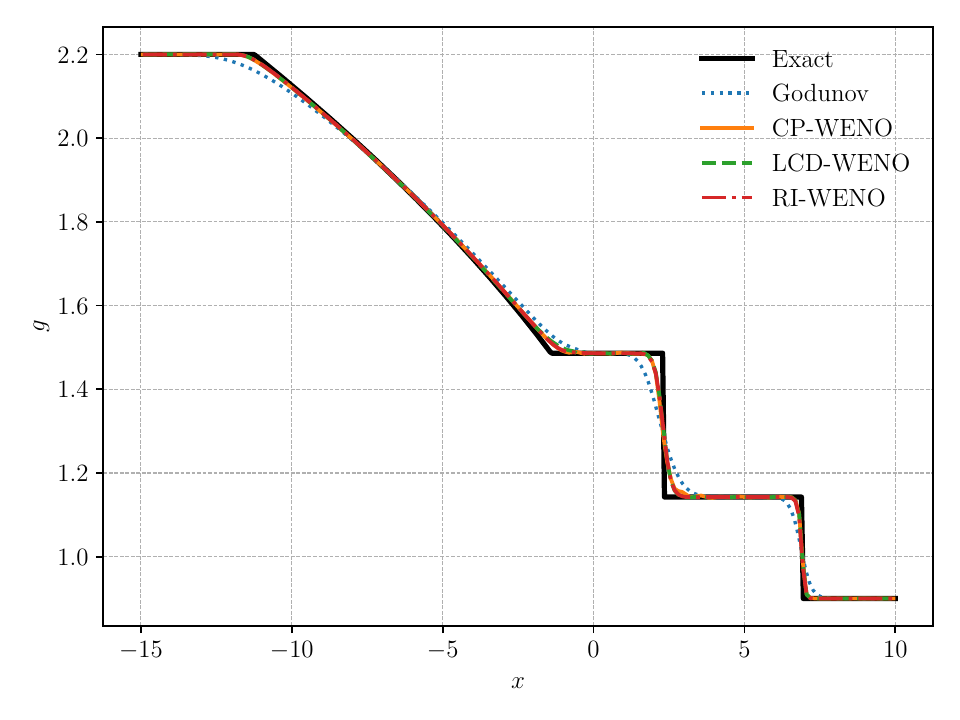}
\caption{Film height ($g$)}
\end{subfigure}

\vspace{0.3cm}

\begin{subfigure}{0.48\textwidth}
\includegraphics[width=\linewidth]{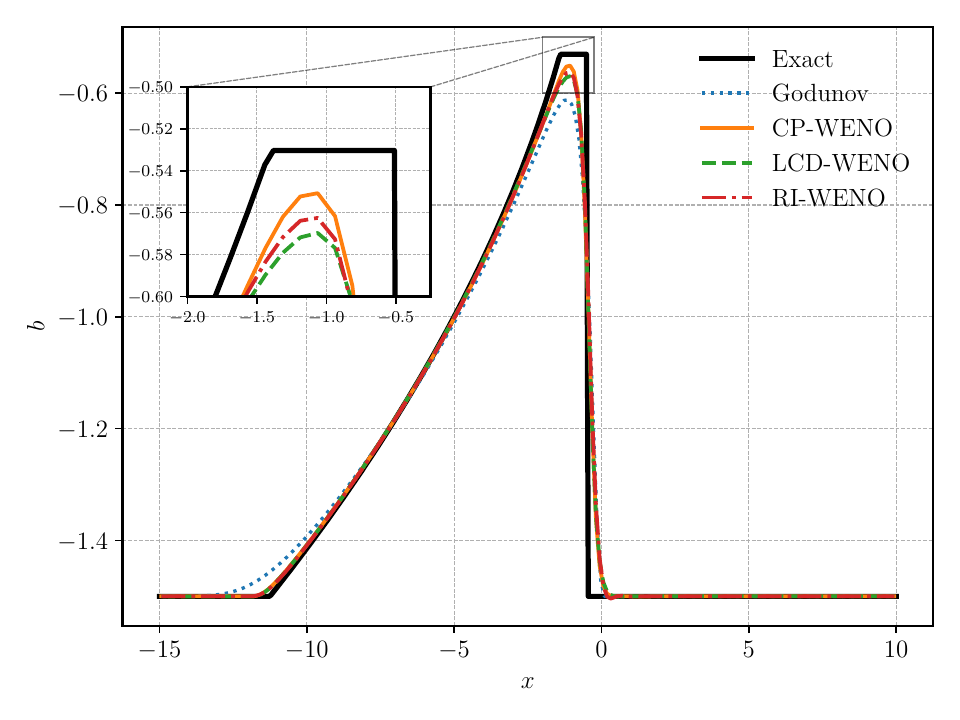}
\caption{Concentration gradient ($b$)}
\end{subfigure}
\hfill
\begin{subfigure}{0.48\textwidth}
\includegraphics[width=\linewidth]{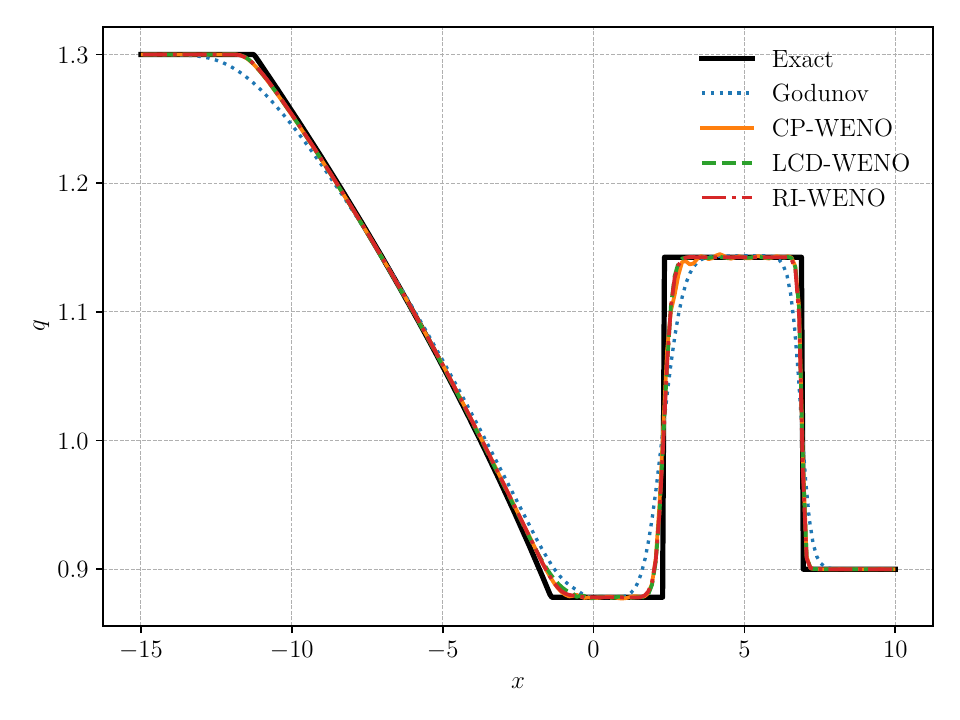}
\caption{Concentration gradient ($q$)}
\end{subfigure}

\caption{Comparison of CP-WENO, LCD-WENO, and RI-WENO with the exact and Godunov solutions for Test Problem~\ref{RP5} at $N_x=200$ and $t=5.0$.}
\label{fig:RP5_all}
\end{figure}

\testproblem{SOD type Riemann Problem \label{RP5}}
This test problem is an SOD-type Riemann problem from \cite{barthwal2026generalized}. The domain is $x \in [-15,10]$ with outflow boundary conditions. The initial profile is given by
$$
(f,b,g,q)(x,0)=
\begin{cases}
(1.0,\,-1.5,\,2.2,\,1.3), & x \le 0, \\
(0.125,\,-1.5,\,0.9,\,0.9), & x > 0.
\end{cases}
$$
The computation is evolved up to the final time $t=5.0$ on a uniform grid with $N_x=200$ cells. The numerical solutions obtained using CP-WENO, LCD-WENO, and RI-WENO are compared against the exact solution and the Godunov scheme in Figure~\ref{fig:RP5_all}. The CP-WENO scheme does not exhibit spurious oscillations in this test case, likely due to the nature of the solution profile. Both LCD-WENO and RI-WENO schemes effectively capture the solution features and closely match the exact solution, with the RI-WENO scheme performing better at the peak of the concentration gradient $b$ compared to the LCD-WENO scheme. Moreover, the RI-WENO scheme turns out to be computationally effective as presented in Table \ref{tab:cpu_time}.
\testproblem{Riemann Problem 3\label{RP3}}

\begin{figure}[htbp]
\centering
\begin{subfigure}{0.48\textwidth}
\includegraphics[width=\linewidth]{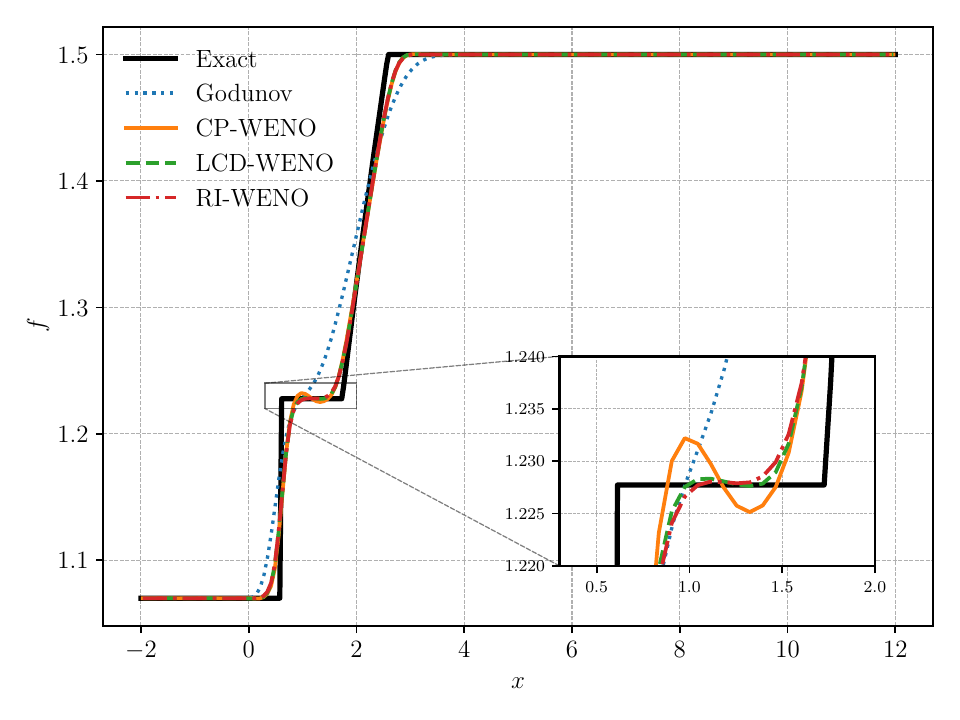}
\caption{Film height ($f$)}
\end{subfigure}
\hfill
\begin{subfigure}{0.48\textwidth}
\includegraphics[width=\linewidth]{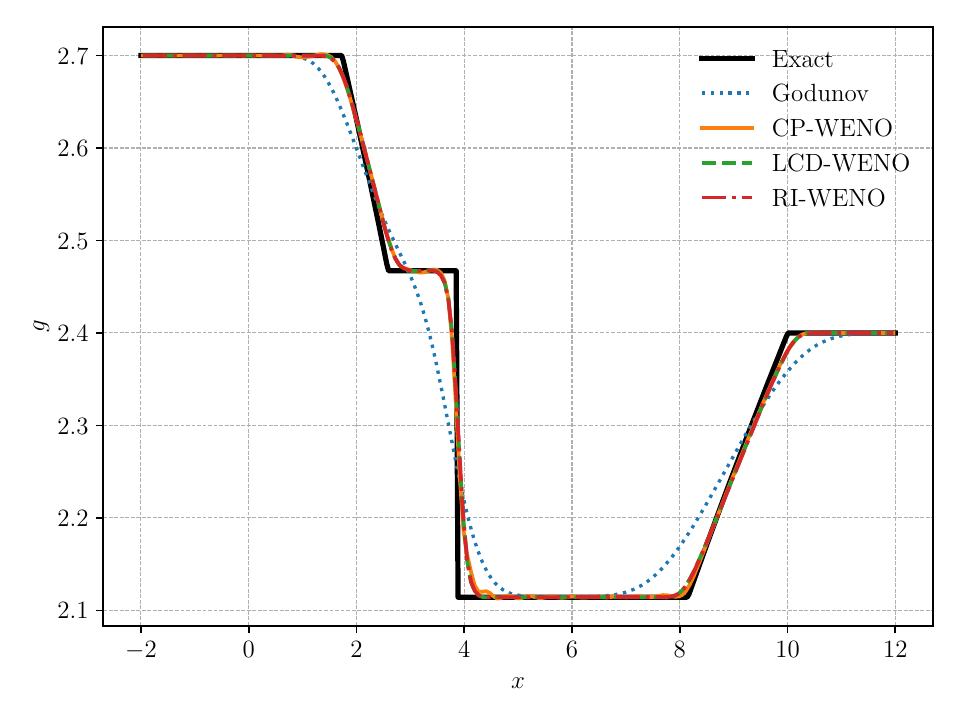}
\caption{Film height ($g$)}
\end{subfigure}

\vspace{0.3cm}

\begin{subfigure}{0.48\textwidth}
\includegraphics[width=\linewidth]{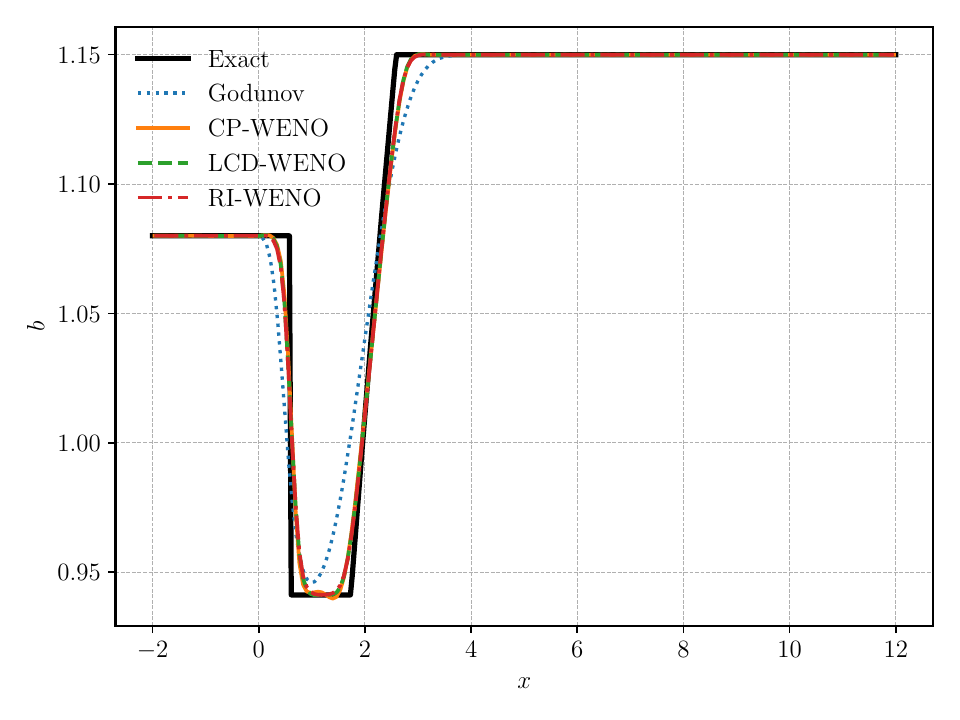}
\caption{Concentration gradient ($b$)}
\end{subfigure}
\hfill
\begin{subfigure}{0.48\textwidth}
\includegraphics[width=\linewidth]{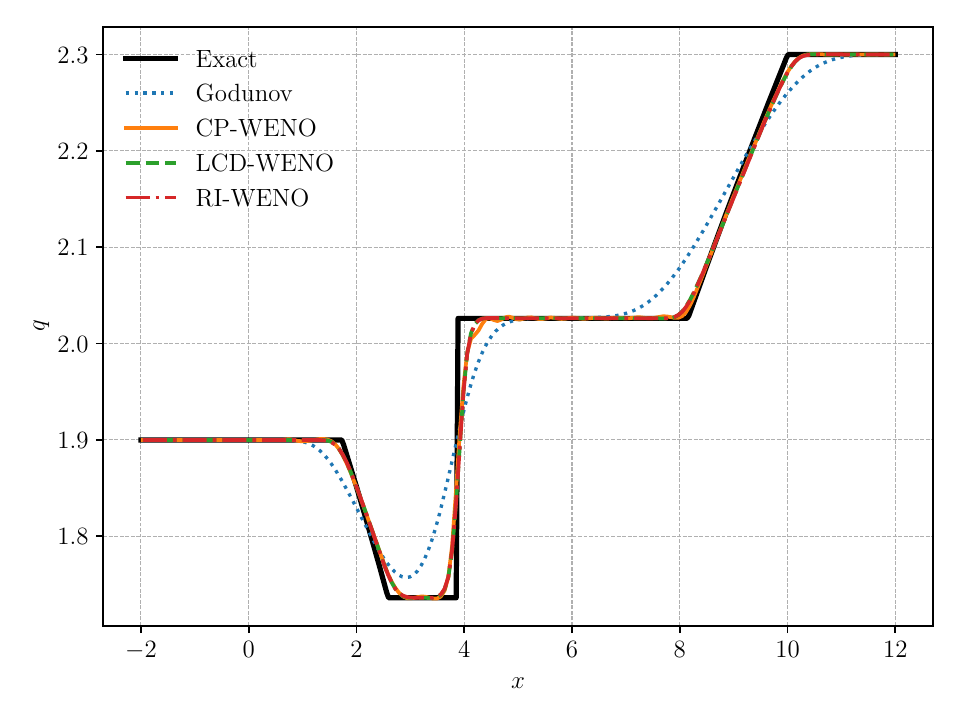}
\caption{Concentration gradient ($q$)}
\end{subfigure}

\caption{Comparison of CP-WENO, LCD-WENO, and RI-WENO with the exact and Godunov solutions for Test Problem~\ref{RP3} at $N_x=200$ and $t=1.0$.}
\label{fig:RP3_all}
\end{figure}

The next test problem is another Riemann problem. The domain is $x \in [-2,12]$ with outflow boundary conditions. We fill the following initial condition,
$$
(f,b,g,q)(x,0)=
\begin{cases}
(1.07,\,1.08,\,2.7,\,1.90), & x \le 0, \\
(1.5,\,1.15,\,2.4,\,2.30), & x > 0.
\end{cases}
$$
The computational domain is discretized using a uniform grid with $N_x=200$ cells, and the solution is evolved up to the final time $t=1.0$. The numerical solutions are compared in Figure~\ref{fig:RP3_all}. We again observe that the CP-WENO scheme exhibits spurious oscillations near the discontinuity, while both LCD-WENO and RI-WENO effectively suppress these oscillations.  Solution profile by RI-WENO scheme almost coincides with that of LCD-WENO while being computationally more efficient as presented in Table~\ref{tab:cpu_time}.

\subsubsection*{Computational cost comparison for Riemann problems}
In this section, we evaluate the speed-up factor of the proposed scheme relative to the CRI-WENO scheme for the Riemann test cases. The speed-up factor is defined as the ratio of the CPU time of a given scheme to that of the CRI-WENO scheme. A speed-up factor less than one indicates that the scheme is faster than CRI-WENO, whereas a value greater than one indicates slower performance.
\begin{table}[htbp]
\centering
\caption{CPU Time Comparison, normalized by CRI-WENO}
\label{tab:cpu_time}
\begin{tabular}{lcccc}
\toprule
\textbf{Problem} & \textbf{CRI-WENO} & \textbf{LCD-WENO} & \textbf{RI-WENO} & \textbf{CP-WENO} \\
\midrule
Test Problem \ref{RP9} & 1.000 & 0.787 & 0.671 & 0.542 \\
Test Problem \ref{RP2} & 1.000 & 0.892 & 0.840 & 0.783 \\
Test Problem \ref{RP5} & 1.000 & 0.892 & 0.853 & 0.797 \\
Test Problem \ref{RP3} & 1.000 & 0.870 & 0.820 & 0.751\\
\bottomrule
\end{tabular}
\end{table}

Table \ref{tab:cpu_time} presents the speed-up factors of various schemes for the Riemann problems. It can be observed that LCD-WENO requires 11–22\% less computational time, RI-WENO requires 15–33\% less time, and CP-WENO requires 20–46\% less time compared to the CRI-WENO scheme.\\\\
\testproblem{Two-Marangoni regime test: Effect of concentration gradients on the film heights}\label{2_marangoni}
To investigate the effect of two-Marangoni numbers on the dynamics of flow, we choose $\alpha=k\beta$ for some $k\in \mathbb{R}$ and choose the initial data of the form
\begin{align*}
f(x,0) &= f_0, 
\qquad c_1(x,0) =
c_{1,-}
+\frac{c_{1,+}-c_{1,-}}{2}
\left(
1+\tanh\!\left(\frac{x-x_0}{\delta}\right)
\right),\\
g(x,0) &= g_0, \qquad
c_2(x,0) =
c_{2,-}
+\frac{c_{2,+}-c_{2,-}}{2}
\left(
1+\tanh\!\left(\frac{x-x_0}{\delta}\right)
\right).
\end{align*}
Physically, this initial condition represents
two initially flat liquid layers with the concentrations $c_1$ and $c_2$ each
take one approximately constant value inside the region $x<x_0$ and another
approximately constant value outside the region $x>x_0$, with a smooth
transition across a narrow region of thickness $\delta$. Thus,
$c_{1,-}$ and $c_{1,+}$ denote the inner and outer concentration
levels for the first layer, while $c_{2,-}$ and $c_{2,+}$ denote the
corresponding levels for the second layer.

Further, we assume that
\[
c_{1,+}>c_{1,-},
\qquad
c_{2,+}>c_{2,-},
\]
so that $c_1$ and $c_2$ increases across the interface $x_0$.

\begin{figure}[htbp]
\centering
\begin{subfigure}{0.49\textwidth}
\includegraphics[width=\linewidth]{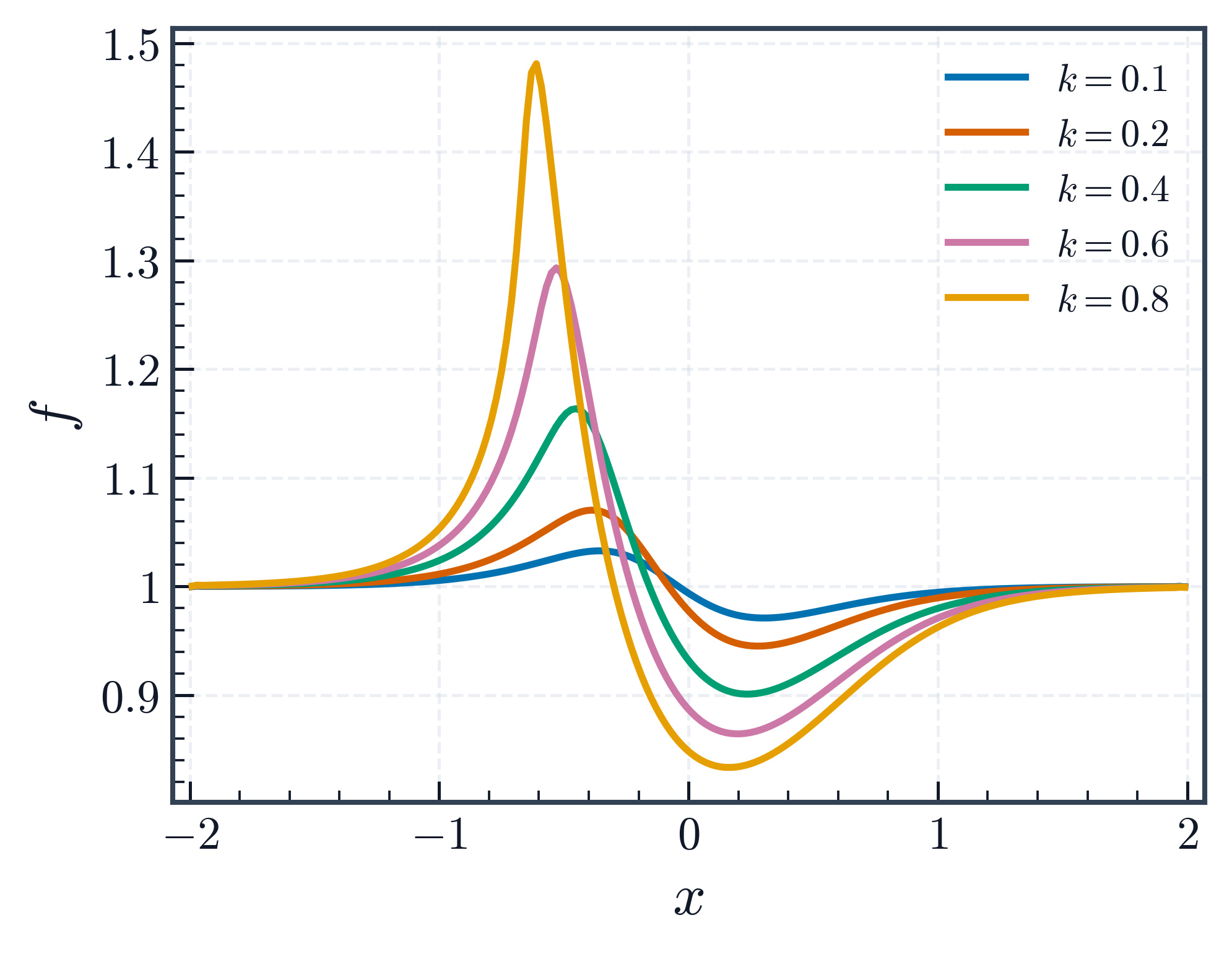}
\caption{Film height ($f$)}
\end{subfigure}\hspace{0pt}
\begin{subfigure}{0.49\textwidth}
\includegraphics[width=\linewidth]{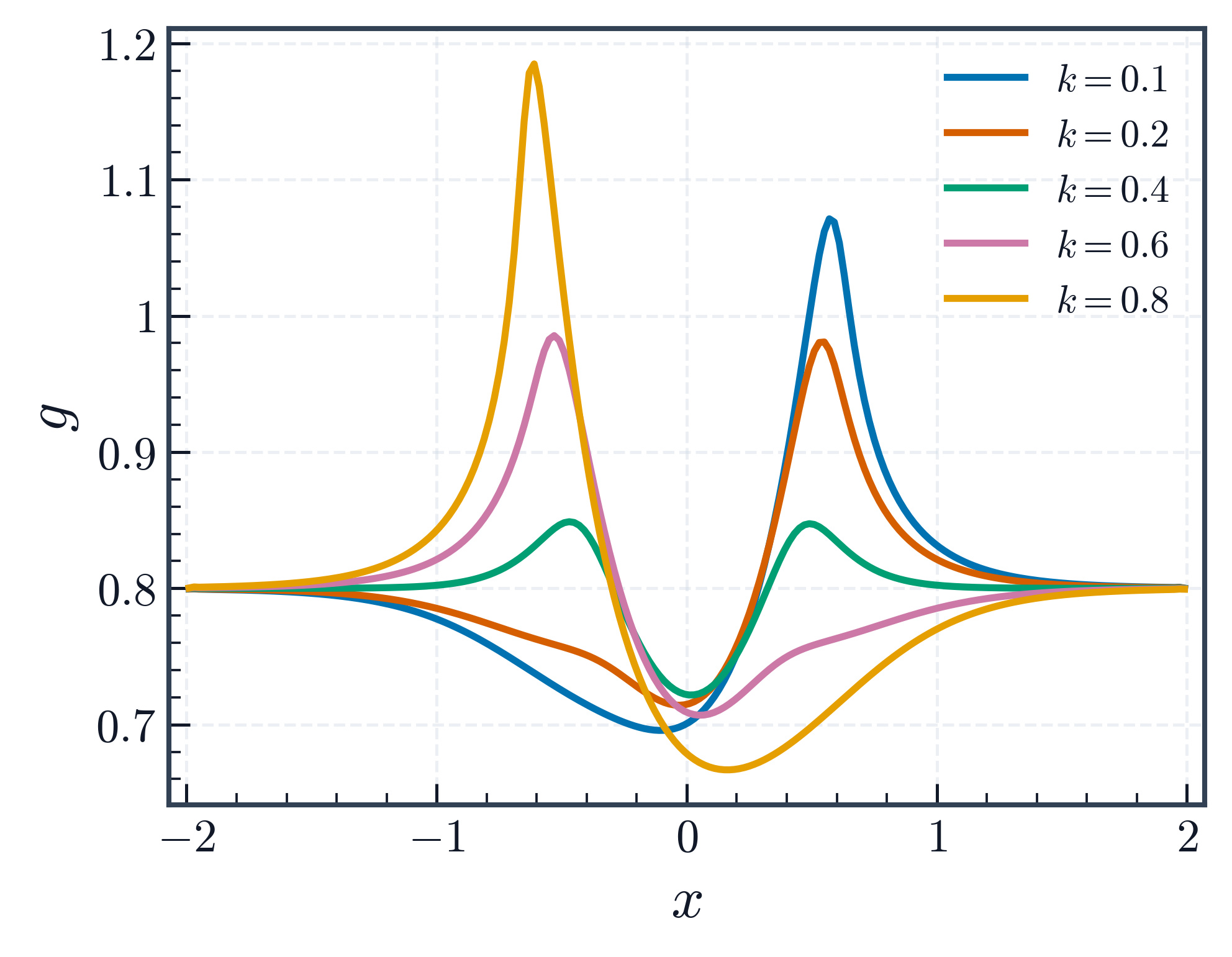}
\caption{Film height ($g$)}
\end{subfigure}

\caption{Comparison of film heights $f$ and $g$ for different values of $k<1$ at $t=1.00$.}
\label{fig:two_marangoni_all_less_k}
\end{figure}

\begin{figure}[htbp]
\centering
\begin{subfigure}{0.49\textwidth}
\includegraphics[width=\linewidth]{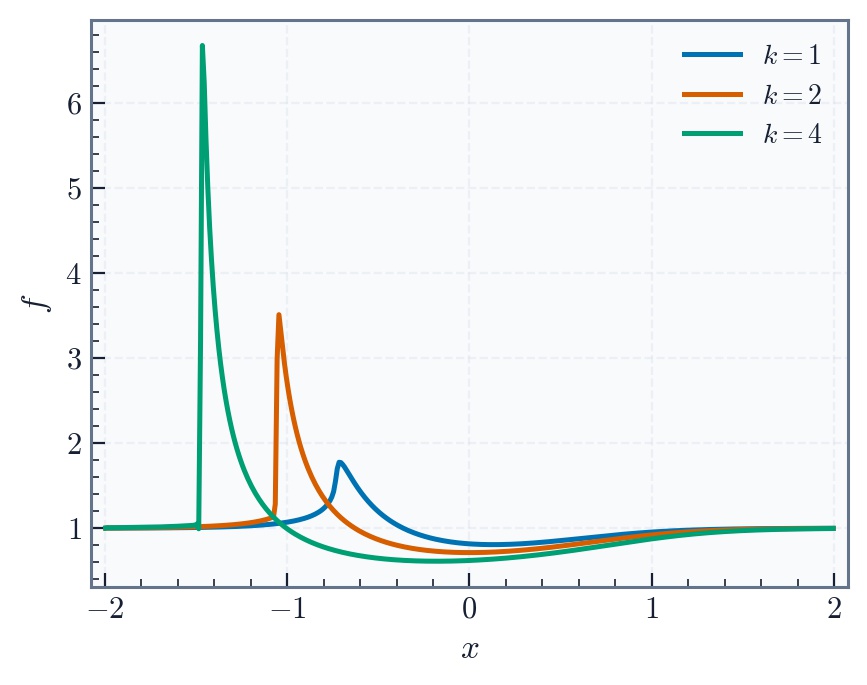}
\caption{Film height ($f$)}
\end{subfigure}\hspace{0pt}
\begin{subfigure}{0.49\textwidth}
\includegraphics[width=\linewidth]{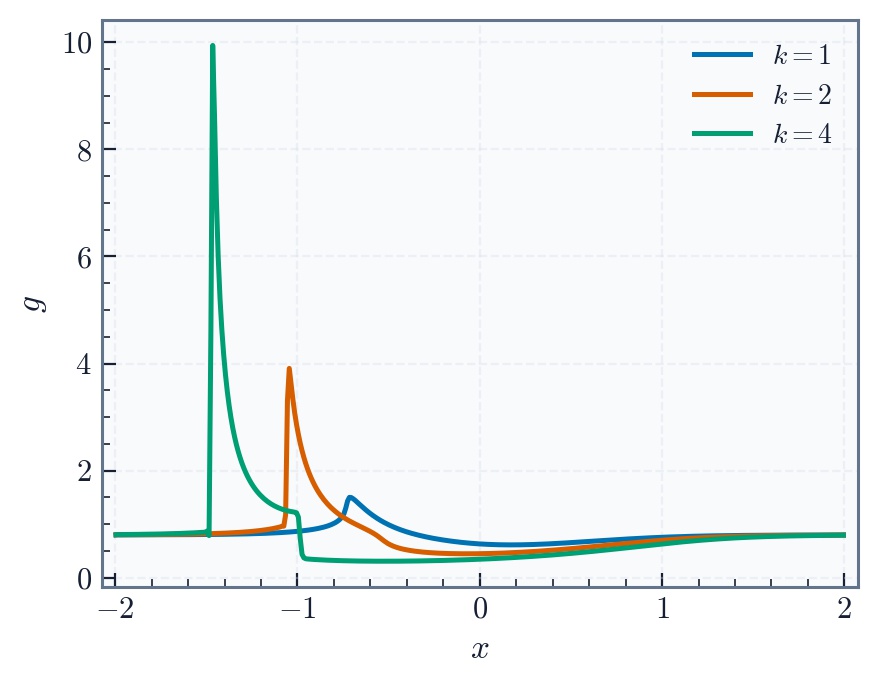}
\caption{Film height ($g$)}
\end{subfigure}
\caption{Comparison of film heights $f$ and $g$ for different values of $k\geq 1$ at $t=1.00$.}
\label{fig:two_marangoni_all_more_k}
\end{figure}

The corresponding initial
data for the conservative variables $(f, b, g, q)$ then are
\begin{equation}
\begin{aligned}
f(x, 0) &= f_0, \qquad
b(x,0) =
-\frac{c_{1,-}-c_{1,+}}{2\delta}\,
\operatorname{sech}^2\!\left(\frac{x-x_0}{\delta}\right),\\
g(x, 0) &= g_0, \qquad
q(x,0) =
\frac{c_{2,+}-c_{2,-}}{2\delta}\,
\operatorname{sech}^2\!\left(\frac{x-x_0}{\delta}\right).
\end{aligned}
\end{equation}
As concrete parameter values, we choose in particular
\begin{align*}
x_0=0, \quad \delta=0.5,\quad f_0=1, \quad g_0=0.8,\quad c_{1,-}=1.0, \quad c_{1,+}=0.6,\quad c_{2,-}=0.2, \quad c_{2,+}=0.6.
\end{align*}
We plot the evolution of the film heights $f$ and $g$ at time $t=1.00$ in Figures~\ref{fig:two_marangoni_all_less_k} and~\ref{fig:two_marangoni_all_more_k} for different values of $k$, distinguishing between the regimes $k<1$ and $k\geq 1$, respectively. For $k<1$, corresponding to $\alpha<\beta$, the Marangoni effect associated with the second layer plays the dominant role in the dynamics. In this regime, the height of the second layer is strongly affected and is transported to the right as $\beta$ increases, with the peak of the post-shock wave governing the flow behaviour. As $k$ increases, however, the Marangoni effect of the first layer becomes dominant. Consequently, the influence of the second Marangoni number is significantly reduced, and the flow is driven increasingly towards the direction of the concentration gradient in the first layer, namely to the left. In addition, the second wave observed in the evolution of the film height $g$ becomes progressively damped as $k$ increases; see Figure~\ref{fig:two_marangoni_all_less_k}. For larger values of $k$, the left-moving shock becomes sharper, and its strength increases considerably. Moreover, a secondary shock appears in the evolution of the film height $g$ for $k=2$, which sharpens for $k=4$; see Figure~\ref{fig:two_marangoni_all_more_k}.

\subsection{Two-dimensional test cases}
As pointed out in section \ref{sec: hyperbolicity}, we are interested in developing solutions for the system \eqref{eq: Main_system_new_reduced} along a normal direction. In particular, one can take the normal vector as $\mathbf{n}=(1/\sqrt{2}, 1/\sqrt{2})$. In other words, the concentration gradient affects are along $y=x$ line. Physically, this implies that the concentration changes at the same rate in the $x$- and $y$-directions, so there is no preferred coordinate direction. Therefore, the steepest variation occurs along the diagonal direction, which reflects a symmetric influence of the concentration gradient in both spatial directions. In this case, the system \eqref{eq: Main_system_new_reduced} can be rewritten as a simplified two-dimensional hyperbolic system of the form
\begin{align}\label{eq: Main_system}
    \dfrac{\partial f}{\partial t}+\dfrac{\partial}{\partial x}\left(\dfrac{\alpha}{2}f^2b\right)+\dfrac{\partial}{\partial y}\left(\dfrac{\alpha}{2}f^2b\right)&=0,\vspace{0.2 cm} \nonumber\\
    \dfrac{\partial b}{\partial t}+\dfrac{\partial}{\partial x}\left(\dfrac{\alpha}{2}fb^2\right)+\dfrac{\partial}{\partial y}\left(\dfrac{\alpha}{2}fb^2\right)&=0,\vspace{0.2 cm} \\
     \dfrac{\partial g}{\partial t}+\dfrac{\partial}{\partial x}\left(\dfrac{\beta g^2q}{2}+\alpha fgb\right)+\dfrac{\partial}{\partial y}\left(\dfrac{\beta g^2q}{2}+\alpha fgb\right)&=0,\vspace{0.2 cm} \nonumber\\
  \dfrac{\partial q}{\partial t}+\dfrac{\partial}{\partial x}\left(\dfrac{\beta gq^2}{2}+\alpha fbq\right)+\dfrac{\partial}{\partial y}\left(\dfrac{\beta gq^2}{2}+\alpha fbq\right)&=0\nonumber
\end{align}
This particular form of the system \eqref{eq: Main_system_new_reduced} allows us to utilize the numerical scheme developed in Section \ref{sec: Numerical_scheme} even for multi-dimensional examples.
\testproblem{Smooth circular perturbation test\label{TWODNEW2}}
In this test case, we consider a smooth circular perturbation problem for the two-dimensional system on the domain $[-4,4]\times[-4,4]$ with periodic boundary conditions in both spatial directions. The simulation is performed up to the final time $t=1.0$. The initial perturbation is centered at $(x_c,y_c)=(-1,-1)$, and the radial coordinate is defined by
$$
r(x,y)=\sqrt{(x-x_c)^2+(y-y_c)^2}.
$$
To construct a smooth interface, we introduce the profile function
$$
\phi(r)=\frac{1}{2}
\left(
1+\tanh\left(2(2-r)\right)
\right).
$$
The initial condition is prescribed as
\begin{equation}
\begin{aligned}
f(x,y,0) &= 1.0, \\
b(x,y,0) &= 0.5 + 0.1\,\phi(r), \\
g(x,y,0) &= 2.0, \\
q(x,y,0) &= 0.5 + 0.2\,\phi(r).
\end{aligned}
\end{equation}
Thus, the variables $b$ and $q$ contain smooth localized perturbations, while the background states of $f$ and $g$ remain spatially constant. The perturbation amplitude in $q$ is chosen larger than that in $b$, leading to different nonlinear transport effects in the two layers.

Since the initial data are smooth and radially symmetric, this test problem is well suited for assessing multidimensional wave propagation. The plots of the primitive variables show that the concentration-gradient profiles generate a Marangoni-driven flow. As a result, the initially constant film heights evolve into nonlinear wave patterns. More precisely, the localized Marangoni forcing redistributes the fluid layers and transports mass throughout the computational domain. Consequently, even though the film heights are initially constant, the solution develops alternating regions of larger and smaller height.

\begin{figure}[htbp]
\centering
\begin{subfigure}{0.48\textwidth}
\includegraphics[width=\linewidth]{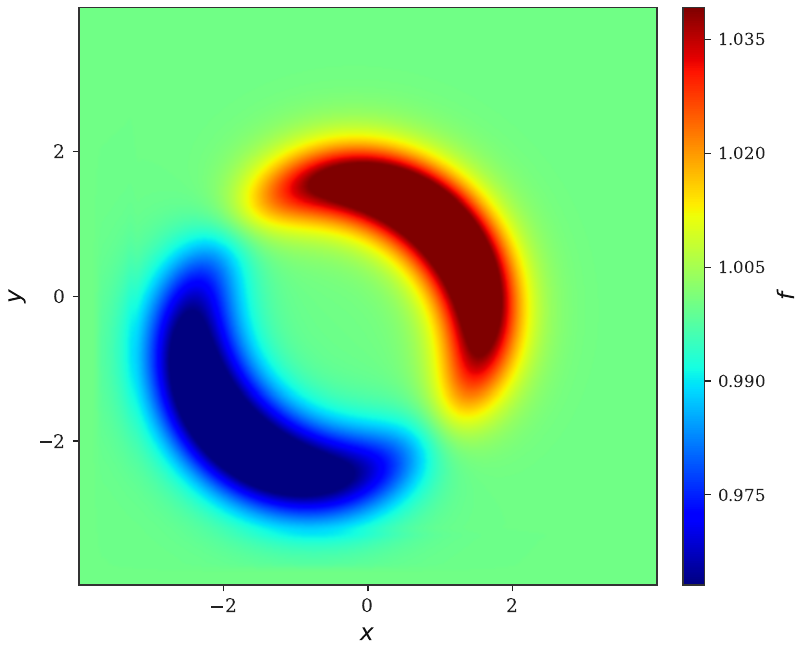}
\caption{Film height ($f$)}
\end{subfigure}
\hfill
\begin{subfigure}{0.48\textwidth}
\includegraphics[width=\linewidth]{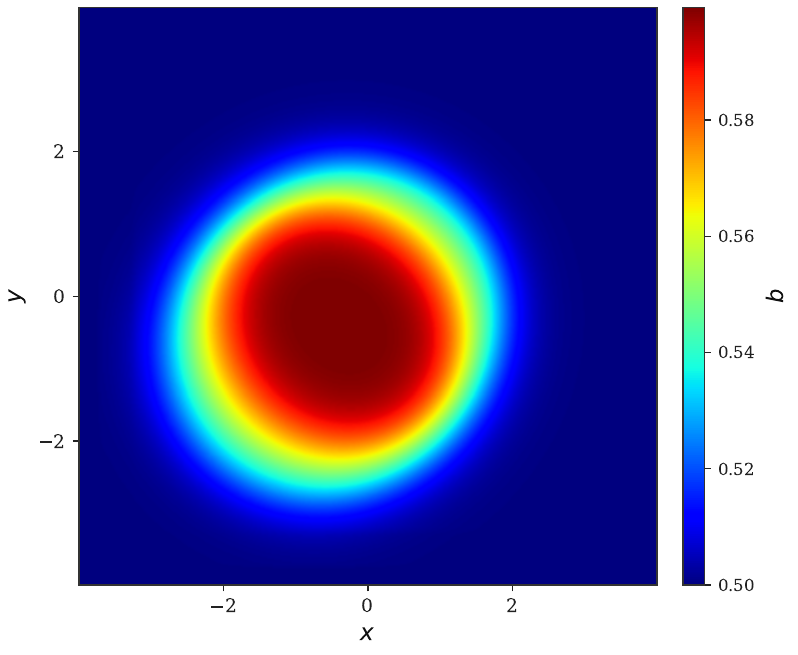}
\caption{Concentration gradient ($b$)}
\end{subfigure}\\
\vspace{0.5cm}

\begin{subfigure}{0.48\textwidth}
\includegraphics[width=\linewidth]{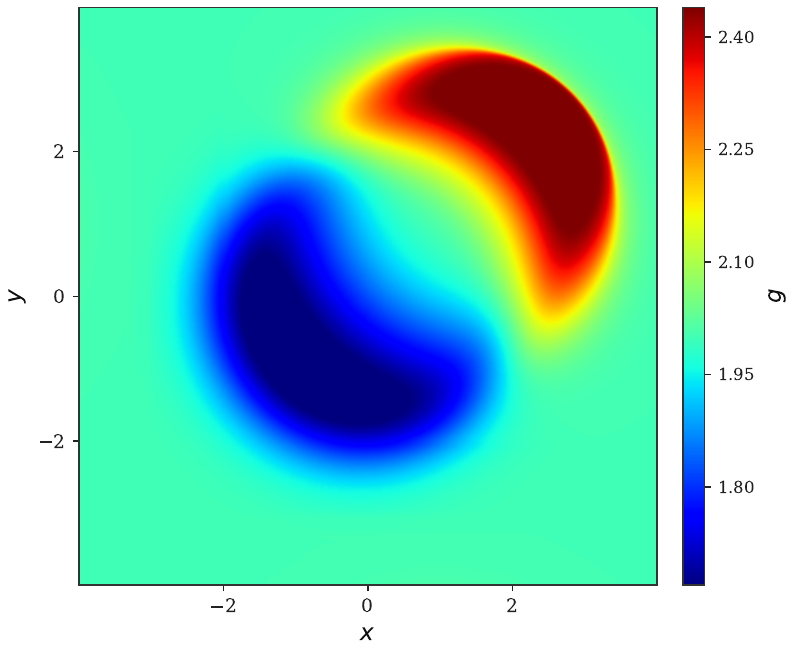}
\caption{Film height ($g$)}
\end{subfigure}
\hfill
\begin{subfigure}{0.48\textwidth}
\includegraphics[width=\linewidth]{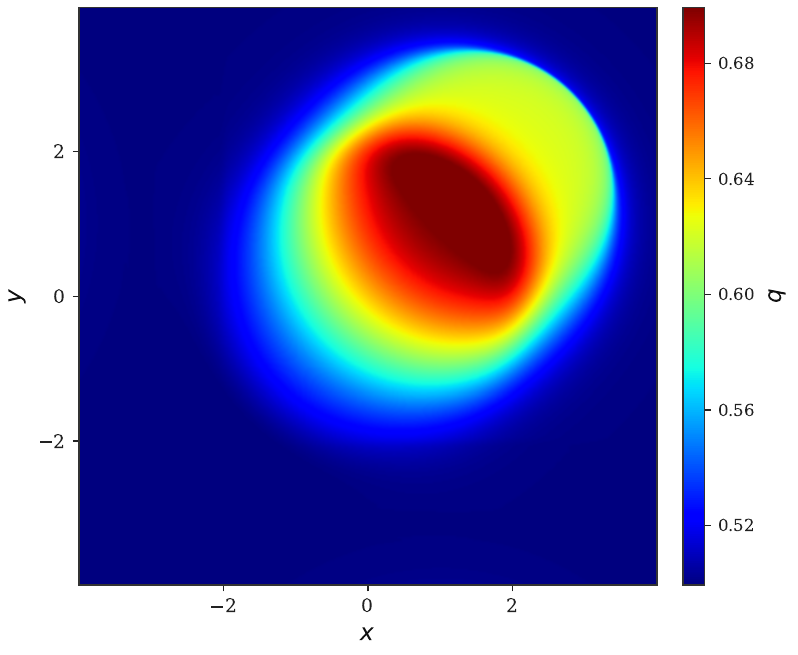}
\caption{Concentration gradient ($q$)}
\end{subfigure}

\caption{Test Problem~\ref{TWODNEW2} with $400\times400$ grids and $t=1.0$.}
\label{fig:TWODNEW2_all}
\end{figure}
\testproblem{Perturbed circular solutal front test\label{TM2D}}
In this test case, we consider the 2D version of the test case considered in Test problem 6 in section \ref{2_marangoni}. In particular, we consider a circular concentration-front initial condition
centered at $(x_0,y_0)$, and define
\[
r(x,y)=\sqrt{(x-x_0)^2+(y-y_0)^2}
\]
and $\theta$ as the polar angle around $(x_0,y_0)$ such that $x-x_0=r\cos \theta$ and $y-y_0=r\sin \theta$. Moreover, for simplicity we assume the same Marangoni scale such that $\alpha=\beta=1$.

The initial data in terms of the film heights $f,g$ and concentrations
$c_1,c_2$ are given by
\begin{align*}
f(x,y,0) &= f_0, 
\qquad c_1(x,y,0) =
c_{1,-}
+\frac{c_{1,+}-c_{1,-}}{2}
\left(
1+\tanh\!\left(\frac{r(x,y)-R(\theta)}{\delta}\right)
\right),\\
g(x,y,0) &= g_0, \qquad
c_2(x,y,0) =
c_{2,-}
+\frac{c_{2,+}-c_{2,-}}{2}
\left(
1+\tanh\!\left(\frac{r(x,y)-R(\theta)}{\delta}\right)
\right).
\end{align*}
Here, $R(\theta)=R_0+\epsilon \cos(2\theta)>0$ for a sufficiently small $\epsilon>0$  denotes the perturbed radius of the circular front and $\delta>0$
is the transition thickness. In this spirit, the circular front actually is a perturbed domain and not a fixed domain.

In contrast to 1D test, we assume that
\[
c_{1,+}<c_{1,-},
\qquad
c_{2,+}>c_{2,-},
\]
so that $c_1$ decreases while $c_2$ increases across the circular
interface. This allows us to consider the competing effects of the concentration difference in each layer.

Interpreting the scalar variables $b$ and $q$ as reduced gradient
amplitudes associated with the circular front, the corresponding initial
data are
\begin{equation}
\begin{aligned}
f(x,y,0) &= f_0, \qquad
b(x,y,0) =
-\frac{c_{1,-}-c_{1,+}}{2\delta}\,
\operatorname{sech}^2\!\left(\frac{r(x,y)-R(\theta)}{\delta}\right),\\
g(x,y,0) &= g_0, \qquad
q(x,y,0) =
\frac{c_{2,+}-c_{2,-}}{2\delta}\,
\operatorname{sech}^2\!\left(\frac{r(x,y)-R(\theta)}{\delta}\right).
\end{aligned}
\end{equation}
As concrete parameter values, we choose in particular
\begin{align*}
(x_0,y_0)=(0,0), \qquad R_0=0.35, \qquad \delta=0.2,\qquad \epsilon=0.03\\
f_0=1, \quad g_0=0.8,\quad c_{1,-}=1.0, \quad c_{1,+}=0.6,\quad c_{2,-}=0.2, \quad c_{2,+}=0.6.
\end{align*}

\begin{figure}[htbp]
\centering
\begin{subfigure}{0.48\textwidth}
\includegraphics[width=\linewidth]{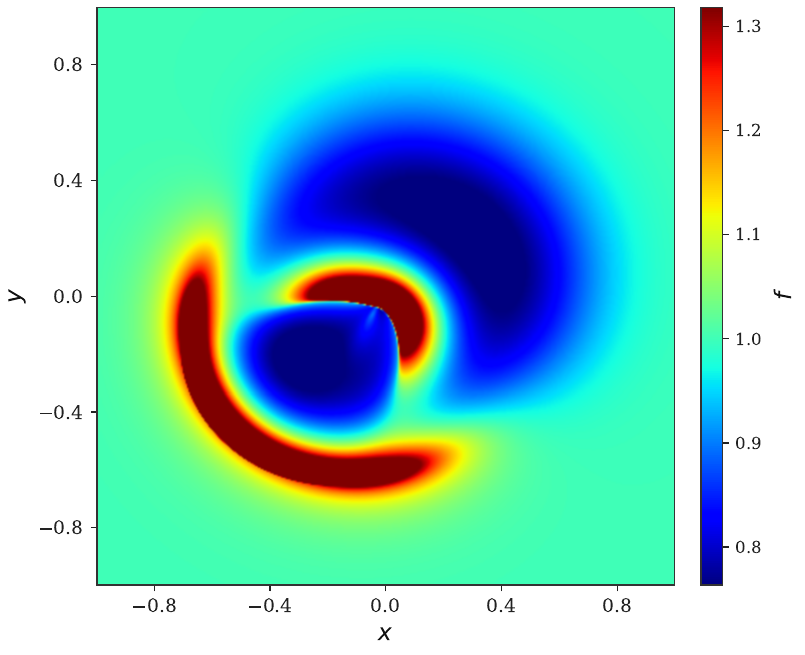}
\caption{Film height ($f$)}
\end{subfigure}
\hfill
\begin{subfigure}{0.48\textwidth}
\includegraphics[width=\linewidth]{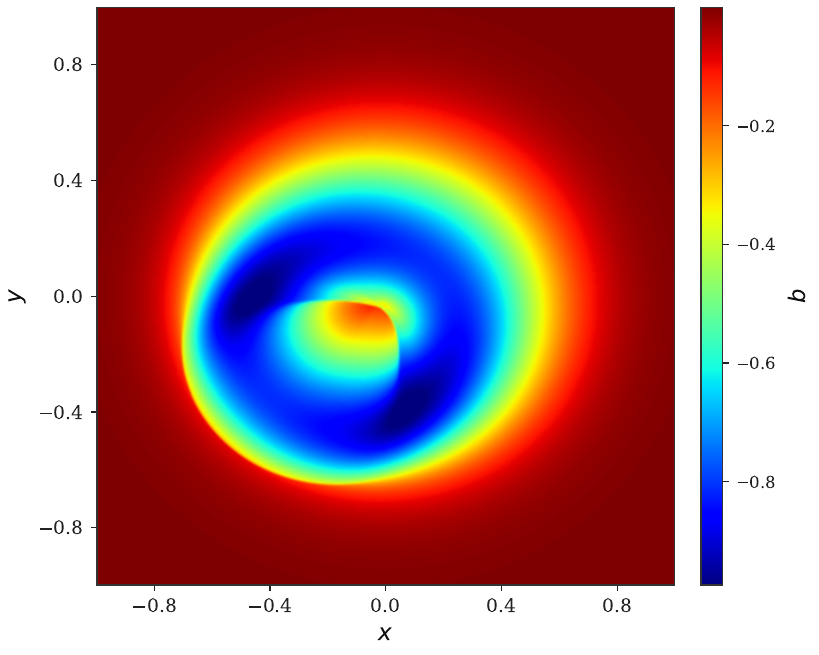}
\caption{Concentration gradient ($b$)}
\end{subfigure}\\
\vspace{0.3cm}

\begin{subfigure}{0.48\textwidth}
\includegraphics[width=\linewidth]{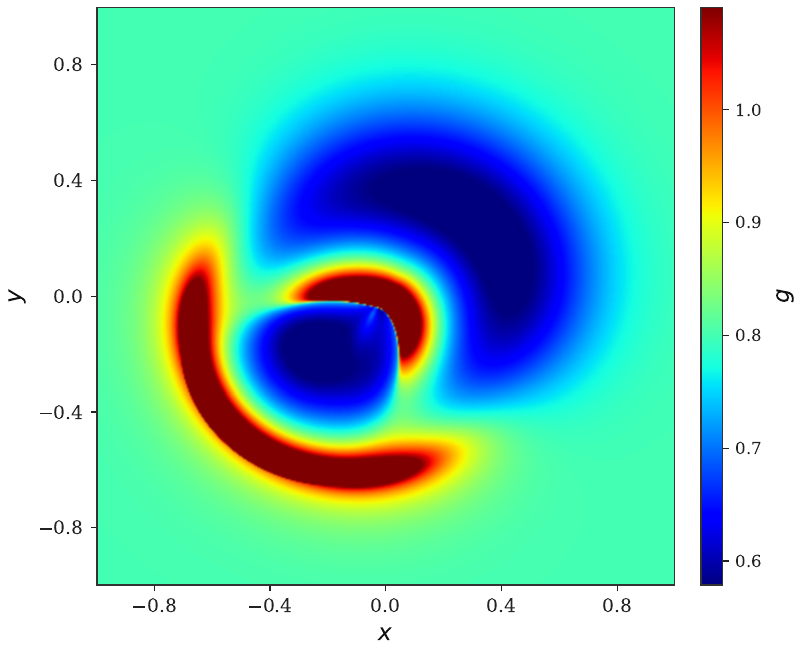}
\caption{Film height ($g$)}
\end{subfigure}
\hfill
\begin{subfigure}{0.48\textwidth}
\includegraphics[width=\linewidth]{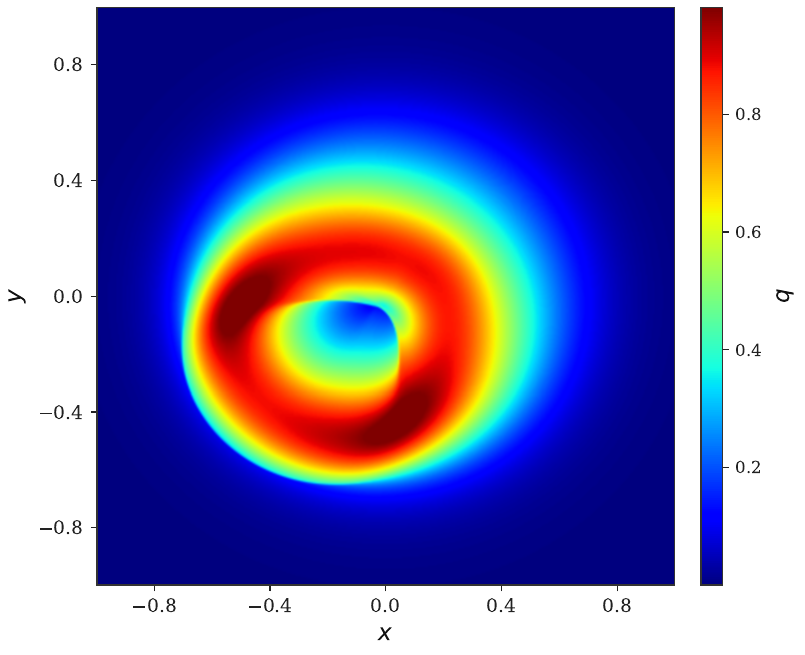}
\caption{Concentration gradient ($q$)}
\end{subfigure}

\caption{Test Problem~\ref{TM2D} with $400\times400$ grids and $t=0.2$.}
\label{fig:TM2D_all}
\end{figure}

We plot the results at time $t=0.2$ for this choice of initial data in Figure \ref{fig:TM2D_all}. From the plots of primitive variables, one can observe that the concentration gradient profiles and induced Marangoni forcing remain concentrated around the thin annular region, while the constant initial film heights evolve as nonlinear waves. In particular, the localized Marangoni forcing redistributes the film heights, and the fluids are thus transported away from some parts of the interface and accumulate in others. That is why, although initial film heights are constant, the solution at time $t=0.2$ develops alternating regions of larger and smaller heights.

\section{Conclusions and future outlook}
In this article, we developed a multidimensional two-layer thin-film model extending the model proposed by Barthwal and Rohde \cite{barthwal2025hyperbolic}. We analyzed its hyperbolicity and Riemann invariant structure. The resulting system admits a Riemann invariant structure in multidimensions along a prescribed spatial direction, which enables the construction of a Riemann invariant-based alternative WENO finite difference scheme. This scheme is computationally cheaper while remaining highly efficient compared with its LCD-based counterparts.

A key aspect of the proposed approach is the careful selection of Riemann invariants so as to reduce the computational cost. Although the system admits a full set of Riemann invariants that can be used as a coordinate system, the corresponding inversion process requires an iterative solver, thereby increasing the computational expense. Therefore, a suitable choice of Riemann invariants, leading to a sparser eigenstructure, is essential for designing more efficient high-order numerical methods.

Since the model possesses an entropy structure and the film heights must remain positive in thin-film flows, a natural next step is to develop a provably positivity-preserving and entropy-stable numerical scheme. Finally, in the present work, we considered only the reduced first-order model and neglected capillarity and diffusivity effects. Another important direction for future work is therefore to analyze the full lubrication model and develop robust numerical schemes for it. \medskip\\\\
\textbf{Acknowledgments} Authors Biswarup Biswas and Rakesh Kumar acknowledge the Mahindra University Super Computing facility for providing computational resources. Biswarup Biswas is supported by the State University Research Excellence (SURE) Scheme of the Anusandhan National Research Foundation (ANRF), India, under File No. SUR/2022/001786. Rakesh Kumar is supported by the Prime Minister Early Career Research Grant (PMECRG) of the ANRF, India, under Grant No. ANRF/ECRG/2025/004846/PMS.
\bibliographystyle{unsrt}
\bibliography{references}

@article {MR4946687,
    AUTHOR = {Arun, K. R. and Kumar, R. and Meena, A. K.},
     TITLE = {Hybrid finite difference {WENO} schemes for the ten-moment
              {G}aussian closure equations with source term},
   JOURNAL = {Wave Motion},
  FJOURNAL = {Wave Motion. An International Journal Reporting Research on
              Wave Phenomena},
    VOLUME = {139},
      YEAR = {2025},
     PAGES = {Paper No. 103614, 22},
      ISSN = {0165-2125,1878-433X},
   MRCLASS = {65M06 (35L65 65M12)},
  MRNUMBER = {4946687},
MRREVIEWER = {Jianzhong\ Chen},
       DOI = {10.1016/j.wavemoti.2025.103614},
       URL = {https://doi.org/10.1016/j.wavemoti.2025.103614},
}

@article {kum_cha_22a,
    AUTHOR = {Kumar, R. and Chandrashekar, P.},
     TITLE = {Multi-level {WENO} schemes with an adaptive
              characteristic-wise reconstruction for system of {E}uler
              equations},
 JOURNAL = {Computers \& Fluids},
    VOLUME = {239},
      YEAR = {2022},
     PAGES = {Paper No. 105386, 24},
      ISSN = {0045-7930},
   MRCLASS = {76M12 (76Nxx)},
  MRNUMBER = {4396537},
       DOI = {10.1016/j.compfluid.2022.105386},
       URL = {https://doi.org/10.1016/j.compfluid.2022.105386},
}

@article {aru_etal_22a,
    AUTHOR = {Arun, K. R. and Dond, A. K. and Kumar, R.},
     TITLE = {{WENO Smoothness Indicator Based Troubled-cell Indicator for Hyperbolic Conservation Laws}},
JOURNAL = {Computers \& Fluids},
    VOLUME = {},
      YEAR = {2023},
     PAGES = {},
      ISSN = {0045-7930},
   MRCLASS = {76M12 (76Nxx)},
  MRNUMBER = {4396537},
       DOI = {},
       URL = {https://doi.org/10.1016/j.compfluid.2022.105386},
}

@article {zhao-etal_20a,
    AUTHOR = {Zhao, G. Y. and Sun, M. B. and Pirozzoli, S.},
     TITLE = {On shock sensors for hybrid compact/{WENO} schemes},
 JOURNAL = {Computers \& Fluids},
    VOLUME = {199},
      YEAR = {2020},
     PAGES = {104439, 17},
      ISSN = {0045-7930},
   MRCLASS = {76M20 (76L05)},
  MRNUMBER = {4057020},
       DOI = {10.1016/j.compfluid.2020.104439},
       URL = {https://doi.org/10.1016/j.compfluid.2020.104439},
}

@article {kum-cha_19a,
    AUTHOR = {Kumar, R. and Chandrashekar, P.},
     TITLE = {Efficient seventh order {WENO} schemes of adaptive order for
              hyperbolic conservation laws},
JOURNAL = {Computers \& Fluids},
    VOLUME = {190},
      YEAR = {2019},
     PAGES = {49--76},
      ISSN = {0045-7930},
   MRCLASS = {65M06 (76M20 76Nxx)},
  MRNUMBER = {3961109},
       DOI = {10.1016/j.compfluid.2019.06.003},
       URL = {https://doi.org/10.1016/j.compfluid.2019.06.003},
}

@article {bal-etal_16a,
    AUTHOR = {Balsara, D. S. and Garain, S. and Shu, C.W.},
     TITLE = {An efficient class of {WENO} schemes with adaptive order},
  JOURNAL = {Journal of Computational Physics},
    VOLUME = {326},
      YEAR = {2016},
     PAGES = {780--804},
      ISSN = {0021-9991},
   MRCLASS = {65M06},
  MRNUMBER = {3554165},
       DOI = {10.1016/j.jcp.2016.09.009},
       URL = {https://doi.org/10.1016/j.jcp.2016.09.009},
}

@article {kum-cha_18a,
    AUTHOR = {Kumar, R. and Chandrashekar, P.},
     TITLE = {Simple smoothness indicator and multi-level adaptive order
              {WENO} scheme for hyperbolic conservation laws},
 JOURNAL = {Journal of Computational Physics},
    VOLUME = {375},
      YEAR = {2018},
     PAGES = {1059--1090},
      ISSN = {0021-9991},
   MRCLASS = {65M06 (76M20)},
  MRNUMBER = {3874573},
       DOI = {10.1016/j.jcp.2018.09.027},
       URL = {https://doi.org/10.1016/j.jcp.2018.09.027},
}

@article{liu1994weighted,
  title={Weighted essentially non-oscillatory schemes},
  author={Liu, X. D. and Osher, S. and Chan, T.},
  journal={Journal of Computational Physics},
  volume={115},
  number={1},
  pages={200--212},
  year={1994}
}

@article{borges2008improved,
  title={An improved weighted essentially non-oscillatory scheme for hyperbolic conservation laws},
  author={Borges, R. and Carmona, M. and Costa, B. and Don, W. S.},
  journal={Journal of Computational Physics},
  volume={227},
  number={6},
  pages={3191--3211},
  year={2008},
  publisher={Elsevier}
}

@article{jiang1996efficient,
  title={{Efficient implementation of weighted ENO schemes}},
  author={Jiang, G. S. and Shu, C.W.},
  journal={Journal of Computational Physics},
  volume={126},
  number={1},
  pages={202--228},
  year={1996},
  publisher={Elsevier}
}

@book{dafermos2005hyperbolic,
 author = {Dafermos, C. M.},
 title = {Hyperbolic conservation laws in continuum physics},
 edition = {4th},
 fseries = {Grundlehren der Mathematischen Wissenschaften},
 series = {Grundlehren der Mathematischen Wissenschaften},
 issn = {0072-7830},
 volume = {325},
 isbn = {978-3-662-49449-3; 978-3-662-49451-6},
 year = {2016},
 publisher = {Berlin: Springer},
 language = {English},
 zbMATH = {6596871},
 Zbl = {1364.35003}
}

@article{barthwal2026generalized,
  title={{A generalized Riemann problem solver for a hyperbolic model of two-layer thin film flow}},
  author={Barthwal, R. and Rohde, C. and Wang, Y.},
  journal={Journal of Scientific Computing},
  volume={106},
  number={1},
  pages={25},
  year={2026},
  publisher={Springer}
}

@article{barthwal2025hyperbolic,
  title={A hyperbolic model for two-layer thin film flow with a perfectly soluble anti-surfactant},
  author={Barthwal, R. and Rohde, C.},
  journal={Accepted for publication at SIAM Journal on Applied Mathematics},
  year={2026}
}

@article{xu2024local,
  title={{Local characteristic decomposition--free high-order finite difference WENO schemes for hyperbolic systems endowed with a coordinate system of Riemann invariants}},
  author={Xu, Z. and Shu, C.W.},
  journal={SIAM Journal on Scientific Computing},
  volume={46},
  number={2},
  pages={A1352--A1372},
  year={2024},
  publisher={SIAM}
}

@article{wu2025finite,
  title={{Finite difference alternative WENO schemes with Riemann invariant-based local characteristic decompositions for compressible Euler equations}},
  author={Wu, Y. and Shu, C.W.},
  journal={Journal of Computational Physics},
  volume={537},
  pages={114104},
  year={2025},
  publisher={Elsevier}
}

@article{barthwal2025existence,
doi = {10.1088/1361-6544/ae4afc},
url = {https://doi.org/10.1088/1361-6544/ae4afc},
year = {2026},
publisher = {IOP Publishing},
volume = {39},
number = {3},
pages = {035006},
author = {Barthwal, R. and Rohde, C. and Sen, A.},
title = {{Existence and stability of the Riemann solutions for a non-symmetric Keyfitz–Kranzer type model}},
journal = {Nonlinearity},
}

@article{pandey2025construction,
  title={{Construction of solutions of the Riemann problem for a two-dimensional Keyfitz-Kranzer type model governing a thin film flow}},
  author={Pandey, A. and Barthwal, R. and Sekhar, T Raja},
  journal={Applied Mathematics and Computation},
  volume={498},
  pages={129378},
  year={2025},
  publisher={Elsevier}
}

@article{jensen1992insoluble,
  title={Insoluble surfactant spreading on a thin viscous film: shock evolution and film rupture},
  author={Jensen, O.E. and Grotberg, J.B.},
  journal={Journal of Fluid Mechanics},
  volume={240},
  pages={259--288},
  year={1992},
  publisher={Cambridge University Press}
}

@book{serre1999systems,
  title={Systems of conservation laws 1: hyperbolicity, entropies, shock waves},
  author={Serre, D.},
  year={1999},
  publisher={Cambridge University Press}
}

@article{craster2009dynamics,
  title={Dynamics and stability of thin liquid films},
  author={Craster, R. V. and Matar, Omar K.},
  journal={Reviews of Modern Physics},
  volume={81},
  number={3},
  pages={1131--1198},
  year={2009},
  publisher={APS}
}

@article{o2002theory,
  title={Theory and modeling of thin film flows},
  author={O’Brien, S. B. G and Schwartz, L. W.},
  journal={Encyclopedia of surface and colloid science},
  volume={1},
  pages={5283--5297},
  year={2002},
  publisher={Marcel Dekker New York}
}

@article{barthwal2023construction,
  title={\relax{Construction of solutions of a two-dimensional Riemann problem for a thin film model of a perfectly soluble antisurfactant solution}},
  author={Barthwal, R. and Raja Sekhar, T. and P. Raja Sekhar, G.},
  journal={Mathematical Methods in the Applied Sciences},
  volume={46},
  number={6},
  pages={7413--7434},
  year={2023},
  publisher={Wiley Online Library}
}

@article{barthwal2022two,
  title={\relax{Two-dimensional non-self-similar Riemann solutions for a thin film model of a perfectly soluble anti-surfactant solution}},
  author={Barthwal, R. and Raja Sekhar, T.},
  journal={Quarterly of Applied Mathematics},
  volume={80},
  number={4},
  pages={717--738},
  year={2022},
  publisher={American Mathematical Society}
}

@article{bertozzi1999undercompressive,
  title={Undercompressive shocks in thin film flows},
  author={Bertozzi, A.L. and M{\"u}nch, A. and Shearer, M.},
  journal={Physica D: Nonlinear Phenomena},
  volume={134},
  number={4},
  pages={431--464},
  year={1999},
  publisher={Elsevier}
}

@article{conn2017simple,
  title={Simple waves and shocks in a thin film of a perfectly soluble anti-surfactant solution},
  author={Conn, J.J.A. and Duffy, B.R. and Pritchard, D. and Wilson, S.K. and Sefiane, K.},
  journal={Journal of Engineering Mathematics},
  volume={107},
  number={1},
  pages={167--178},
  year={2017},
  publisher={Springer}
}

@article{conn2016fluid,
  title={Fluid-dynamical model for antisurfactants},
  author={Conn, J.J.A. and Duffy, B. R. and Pritchard, D. and Wilson, Stephen K. and Halling, P. J. and Sefiane, K.},
  journal={Physical Review E},
  volume={93},
  number={4},
  pages={043121},
  year={2016},
  publisher={APS}
}

@article{levy2006motion,
  title={The motion of a thin liquid film driven by surfactant and gravity},
  author={Levy, R. and Shearer, M.},
  journal={SIAM Journal on Applied Mathematics},
  volume={66},
  number={5},
  pages={1588--1609},
  year={2006},
  publisher={SIAM}
}

@article{cook2008shock,
  title={Shock solutions for particle-laden thin films},
  author={Cook, B.P. and Bertozzi, A. L. and Hosoi, A.E.},
  journal={SIAM Journal on Applied Mathematics},
  volume={68},
  number={3},
  pages={760--783},
  year={2008},
  publisher={SIAM}
}

@article{myers1998thin,
  title={Thin films with high surface tension},
  author={Myers, T.G.},
  journal={SIAM Review},
  volume={40},
  number={3},
  pages={441--462},
  year={1998},
  publisher={SIAM}
}

@article{matar2004rupture,
  title={Rupture of a surfactant-covered thin liquid film on a flexible wall},
  author={Matar, O.K. and Kumar, S.},
  journal={SIAM Journal on Applied Mathematics},
  volume={64},
  number={6},
  pages={2144--2166},
  year={2004},
  publisher={SIAM}
}

@article{jiang2013alternative,
  title={{An alternative formulation of finite difference weighted ENO schemes with Lax--Wendroff time discretization for conservation laws}},
  author={Jiang, Y. and Shu, C.W. and Zhang, M.},
  journal={SIAM Journal on Scientific Computing},
  volume={35},
  number={2},
  pages={A1137--A1160},
  year={2013},
  publisher={SIAM}
}

@article{balsara2025efficient,
  title={{Efficient alternative finite difference WENO schemes for hyperbolic conservation laws}},
  author={Balsara, D. S. and Bhoriya, D. and Shu, C. W. and Kumar, H.},
  journal={Communications on Applied Mathematics and Computation},
  volume={7},
  number={6},
  pages={2189--2242},
  year={2025},
  publisher={Springer}
}

@article{gottlieb2001strong,
  title={{Strong stability-preserving high-order time discretization methods}},
  author={Gottlieb, S. and Shu, C.W. and Tadmor, E.},
  journal={SIAM Review},
  volume={43},
  number={1},
  pages={89--112},
  year={2001},
  publisher={SIAM}
}
\end{document}